\documentclass{amsart}
\usepackage{changepage}
\usepackage{enumerate} 
\usepackage{graphicx}
\usepackage{hyperref}  
\usepackage{tikz}
\usepackage{tkz-graph}
\usetikzlibrary{matrix,arrows,decorations.markings}
\tikzset{->-/.style={decoration={markings,mark=at position #1 with {\color{black}\arrow{>}}},postaction={decorate,very thick}}}
\tikzstyle{vertex}=[circle, draw, inner sep=0pt, minimum size=6pt]
\newcommand{\vertex}{\node[vertex]}
\pgfrealjobname{ongraphsofHeckeOperators}

\newtheorem*{theorem*}{Theorem}
\newtheorem{theorem}{Theorem}[section]
\newtheorem{lemma}[theorem]{Lemma}
\newtheorem{proposition}[theorem]{Proposition}
\newtheorem{corollary}[theorem]{Corollary}
\theoremstyle{definition}
\newtheorem{definition}[theorem]{Definition}
\newtheorem{example}[theorem]{Example}

\newtheorem{remark}[theorem]{Remark}
\newcommand{\Q}{\mathbb{Q}}
\newcommand{\Z}{\mathbb{Z}}

\allowdisplaybreaks
\usepackage{comment}

\newcommand{\C}{\mathbb{C}}
\newcommand{\norm}[1]{|\!| #1 |\!|}

\newcommand{\E}{\mathcal{E}}
\newcommand{\F}{\mathcal{F}}
\newcommand{\Line}{\mathcal{L}}

\newcommand{\Fq}{\mathbb{F}_q}
\numberwithin{equation}{section}
\usepackage[matrix,arrow]{xy}
\begin{document}

\raggedbottom

\title{On graphs of Hecke operators}


\author{Roberto  Alvarenga}
\address{Estrada Dona Castoria, IMPA, Rio de Janeiro-Brazil}
\curraddr{}
\email{robertoa@impa.br}





\maketitle

\noindent
\textbf{Abstract.} The graph of a Hecke operator encodes all information about the action of this operator on automorphic forms. Let $X$ be a curve over $\mathbb{F}_q$, $F$ its function field and $\mathbb{A}$ the adele ring of $F$. 
In this paper we will exhibit the first properties for the graph of Hecke operators for $\mathrm{GL}_n(\mathbb{A}),$ for every $n \geq 1.$  This includes a description of the graph in terms of coherent sheaves on $X.$ We provide a numerical condition for two vertices to be connected by an edge. Moreover, we describe how to calculate  these graphs in the case of the projective line $X = \mathbb{P}^1(\mathbb{F}_q).$ \\

\tableofcontents

\section*{Introduction}


Zagier observes in \cite{zagier} that if the kernel of certain operators on automorphic forms turns out to be an unitarizable representation, a formula of Hecke implies the Riemann hypothesis. Zagier calls the elements of this kernel \textit{toroidal automorphic forms} (see \cite{zagier} or \cite{oliver-toroidal} for a precise definition). Moreover, Zagier asks what happens if $\Q$ is replaced by a global function field and remarks that  the space of unramified toroidal automorphic forms can be expected to be finite dimensional. 

Motivated by Zagier's question, Lorscheid proves in his Ph.D. thesis, among other things, the following :

$\bullet$(\cite{oliver-graphs} Theorem 10.9) The space of unramified toroidal automorphic forms for a global function field is finite dimensional.

$\bullet$(\cite{oliver-toroidal} Theorem 7.7) There are no nontrivial unramified toroidal automorphic forms for rational function fields. 

$\bullet$(\cite{oliver-elliptic} Theorem 7.12) Let $F$ be the function field of an elliptic curve over a finite field with $q$ elements and class number $h$, let $s+\frac{1}{2}$ be a zero of the zeta function of $F$. If the characteristic is not $2$ or $h\neq q+1$, the space of unramified toroidal automorphic forms is one dimensional and spanned by the Eisenstein series of weight $s$. 
 
The main tool for the proofs of the above theorems is the theory of graphs of Hecke operators, cf. \cite{oliver-graphs}. In this paper, our goal is to continue the study of graphs of Hecke operators and to generalize some results of Lorscheid in \cite{oliver-graphs} from $\mathrm{PGL}_2$ to $\mathrm{GL}_n$. We summarize our main results in the following.

Let $F$ be a function field over $\mathbb{F}_q$, $\mathbb{A}$ its adele ring and $K = \mathrm{GL}_n(\mathcal{O}_{\mathbb{A}}).$ For any right $K$-invariant Hecke operator $\Phi$, there are $m_1, \ldots,m_r \in \C^{*}$ and $[g_1], \ldots, [g_r] \in \mathrm{GL}_n(F) \setminus \mathrm{GL}_n(\mathbb{A}) / K $ such that for all automorphic forms $f$
$$\Phi(f)(g) = \sum_{i=1}^{r} m_i f(g_i),$$
see \ref{propfund}. 
From the above notation we say that there is an edge from $[g]$ to $[g_i]$ with multiplicity $m_i$. If 
$$\Phi = \Phi_{x,r} := \mathrm{char} \left( K \left( \begin{array}{cc} \pi_x I_r & 0  \\
0 & I_{n-r} 
\end{array} \right) K \right), $$
where $x$ is a place on $F$ and $\pi_x$ its uniformizer, we call each $[g_i]$ the $\Phi_{x,r}-$neighbors of $[g]$; see sections \ref{section1} and \ref{section2} for precise definitions. Let $\kappa(x)$ be the residue field of the place $x$. Our first main result is the following.  

\begin{theorem*}[Theorem \ref{maintheorem}] The $\Phi_{x,r}-$neighbors of $[g] \in G(F)\setminus G(\mathbb{A}) / K$ are the classes $[g \xi_w]$ where $\xi_w$ corresponds to elements  $w \in \mathrm{Gr}(n-r,n)(\kappa(x))$. The multiplicity of an edge of  $[g]$ to $[g']$ equals the number of $w \in \mathrm{Gr}(n-r,n)(\kappa(x))$ such that $[g \xi_w] = [g'].$ The multiplicities of the edges originating in $[g]$ sum up to $\# \mathrm{Gr}(n-r,n)(\kappa(x)).$
\end{theorem*}

By a theorem due to Weil, there is a bijection of $ G(F)\setminus G(\mathbb{A}) / K$ with the set of isomorphism classes of rank-$n$ vector bundles on the smooth, projective and geometrically irreducible curve $X$ whose function field is $F$. This theorem allows us to determine the action of an unramified Hecke operator in terms of the equivalence classes of exact sequences of coherent sheaves on $X.$ Namely, we consider exact sequences of the form
$$0 \longrightarrow \E' \longrightarrow \E \longrightarrow \mathcal{K}_{x}^{\oplus r} \longrightarrow 0$$  
where $\E',\E$ are vector bundles, $x$ is a closed point of $X,$ and $\mathcal{K}_{x}^{\oplus r}$ is the skyscraper sheaf on $x$ whose stalk is $\kappa(x)^{\oplus r}$. Let $m_{x,r}(\E,\E')$ be the number of isomorphism classes of exact sequences 
$$0 \longrightarrow \E'' \longrightarrow \E \longrightarrow \mathcal{K}_{x}^{\oplus r} \longrightarrow 0$$
with fixed $\E$ such that $\E'' \cong \E'.$ We denote by  $\mathcal{V}_{x,r}(\E)$ the set of $\big(\E,\E',m_{x,r}(\E,\E')\big)$ such that exists an exact sequence of the type as  above, i.e.\  $m_{x,r}(\E,\E') \neq 0$.   For a subbundle $\mathcal{E'}$ of a bundle $\E$, we define 
$$\delta(\mathcal{E}', \mathcal{E}) := \mathrm{rk}(\E) \deg(\mathcal{E}') - \mathrm{rk}(\E') \deg(\mathcal{E})$$
and
$$\delta_k(\mathcal{E}) := \sup_{\substack{\mathcal{E'} \hookrightarrow \mathcal{E} \\ k-subbundle} } \delta(\mathcal{E}', \mathcal{E}),$$
for $k=1, \ldots, n-1,$ where $n = \mathrm{rk}(\E).$ For precise definitions and notations, see sections \ref{section3} and \ref{section4}. 
 
\begin{theorem*}[Theorem \ref{2maintheorem}] Let $\E$ be a $n$-bundle on $X$ and $x$ be a closed point at $X$ of degree $|x|$. If $m_{x,r}(\E,\E') \neq 0$, then 
$$\delta_k(\E') \in \big\{ \delta_{k}(\E) -k |x|(n-r), \delta_{k}(\E) -k |x|(n-r) + n, \ldots,  \delta_{k}(\E)+ r  |x|k\big\}$$
for every $k=1, \ldots, n-1$. 
\end{theorem*}

In section \ref{section5}, we explain an algorithmic way to calculate the graphs of unramified Hecke operators for a rational function field.

All the results in the first four sections are available for any smooth projective and geometrically irreducible curve over a finite field. Only in the last section do we specialize to the projective line. 

The paper is organized as follows. In the first section, we introduce our notation and provide the basic definitions. In particular, we extend Lorscheid's definition of the graph of a Hecke operator from $\mathrm{PGL}_2$ to $\mathrm{GL}_n.$ In section \ref{section2}, we specialize the theory for unramified Hecke operators $\Phi_{x,r}$ and prove Theorem \ref{maintheorem}. In section \ref{section3}, we describe these graphs in geometric terms, i.e.\ we use the structure of the curve given by $F$ to describe the graph in terms of isomorphism class of short exact sequences of coherent sheaves over this curve. In  section \ref{section4}, we define a numerical invariant for a vector bundle and prove a necessary condition for two vector bundles to be connected by an edge. In section \ref{section5}, we describe an algorithm to calculate the graphs $\Phi_{x,r}$ of Hecke operators for every place $x$, every $n$ and every $r$ when the curve is the projective line.

\section{Background} \label{section1}

In this first section, we set up the notation that is used throughout the paper and introduce graphs of Hecke operators and their first properties.

Let $F$ be a global function field over a finite field $\mathbb{F}_q$ where $q$ is a prime power, i.e. the function field of a geometrically irreducible smooth projective curve $X$ over $\Fq$. Let $g$ be the genus of $X.$ Let $|X|$ the set of closed points of $X$ or, equivalently, the set of places in $F$. For $x \in |X|$, we denote $F_x$ the completion of $F$ at $x$, by $\mathcal{O}_x$ its integers, by $\pi_x \in \mathcal{O}_x$ (we can suppose $\pi_x \in F$) a uniformizer and by $q_x$ the cardinality of the residue field $\kappa(x):=\mathcal{O}_x/(\pi_x) \cong \mathbb{F}_{q_x}.$ Let $|x|$ be the degree of $x$ which is defined by the extension field degree $[\kappa(x) : \Fq]$, in other words $q_x = q^{|x|}$. Let $| \cdot |_x$ the absolute value of $F_x$ (resp. $F$) such that $|\pi_x|_x = q_{x}^{-1}.$

Let $\mathbb{A}$ be the adele ring of $F$ and $\mathbb{A}^{*} $ the idele group. Put $\mathcal{O}_{\mathbb{A}} = \prod \mathcal{O}_x $ where the product is
taken over all places $x$ of $F$. The idele norm is the quasi-character $|\cdot | : \mathbb{A}^{*} \rightarrow  \C^{*}$ that sends
an idele $(a_x) \in \mathbb{A}^{*}$ to the product $\prod |a_x |_x$ over all local norms. By the product formula, this defines a quasi-character on the idele class group $\mathbb{A}^{*}/ F^{*}.$
We think of $F_x$ being embedded into the adele ring $\mathbb{A}$ by sending an element $a \in F_x$ to the adele $(a_y)_{y \in |X|}$ with $a_x = a$ and $a_y = 0$ for $y \neq x$. Not quite compatible with this embedding, we think of the unit group $F_{x}^{*}$ as a subgroup of the idele group $\mathbb{A}^{*}$  by sending an element $b$ of $F_{x}^{*}$
to the idele $(b_y)$ with $b_x = b$ and $b_y = 1$ for $y \neq x$. We will explain in case of ambiguity, which of these embeddings we use.

Let $G(\mathbb{A}):= \mathrm{GL}_n(\mathbb{A}),$ $Z(\mathbb{A})$ be the center of $G(\mathbb{A})$, $G(F):= \mathrm{GL}_n(F)$ and $K:= \mathrm{GL}_n(\mathcal{O}_{\mathbb{A}})$ the standard maximal compact open subgroup of $G(\mathbb{A})$. Note that $G(\mathbb{A})$ comes together with an adelic topology that turns $G(\mathbb{A})$ into a locally compact group. We fix the Haar measure on $G(\mathbb{A})$ for which $\mathrm{vol}(K)=1.$ The topology of $G(\mathbb{A})$ has a neighborhood basis $\mathcal{V}$ of the identity matrix that is given by all subgroups
$$K' = \prod_{x \in |X|} K_{x}' < \prod_{x \in |X|}K_x = K$$
where $K_x := \mathrm{GL}_n(\mathcal{O}_x)$, such that for all $x \in |X|$ the subgroup $K_{x}'$ of $K_x$ is open and consequently of finite index and such that $K_{x}^{'}$ differs from $K_x$ only for a finite number of places. 

Consider the space $C^0(G(\mathbb{A}))$ of continuous functions $f: G(\mathbb{A}) \rightarrow \C.$ Such a function is called \textit{smooth} if it is locally constant. $G(\mathbb{A})$ acts on $C^0(G(\mathbb{A}))$ through the \textit{right regular representation} $\rho : G(\mathbb{A}) \rightarrow \mathrm{Aut}(C^0(G(\mathbb{A})))$, that is defined by right translation of the argument: $(g.f)(h):= (\rho(g)f)(h) := f(hg)$ for $g,h \in G(\mathbb{A})$ and $f \in C^0(G(\mathbb{A})).$

A function $f \in C^{0}(G(\mathbb{A}))$ is called $K-$\textit{finite} if the complex vector space that is generated by $\{k.f\}_{k \in K}$ is finite dimensional. 

Let $H$ be a subgroup of $G(\mathbb{A})$. We say that $f \in C^{0}(G(\mathbb{A}))$ is \textit{left} or \textit{right} $H-$\textit{invariant} if for all $h \in H$ and $g \in G(\mathbb{A})$, 
$ f(hg) = f(g) \text{ or } f(gh)= f(g)$,  respectively.
If $f$ is right and left $H-$invariant, it is called \textit{bi}-$H-$\textit{invariant}.

We embed $G(\mathbb{A}) \hookrightarrow \mathbb{A}^{n^2 +1}$ via $g \mapsto (g, \det(g)^{-1}).$ We define a local height $\norm{ g_x}_x$ on $G(F_x) :=\mathrm{GL}_n(F_x)$ by 
restricting the height function 
$$(v_1, \ldots, v_{n^2+1}) \mapsto \mathrm{max}\{|v_1|_x, \ldots, |v_{n^2+1}|_x\}$$ 
on $F_{x}^{n^2+1}$. We note that $\norm{ g_x }_x \geq 1$ and that $\norm{ g_x}_x =1$ if $g_x \in K_x.$ We define the global height $\norm{ g}$ to be the product of the local heights. We say that $f \in C^{0}(G(\mathbb{A}))$ is of \textit{moderate growth} if there exists constants $C$ and $N$ such that
$$|f(g)|_{\C} \leq C \norm{ g}^{N}$$  
for all $g \in G(\mathbb{A}).$


\begin{definition} The complex vector space $\mathcal{H}$ of all smooth compactly supported functions $\Phi : G(\mathbb{A}) \rightarrow \C$ together with the convolution product
$$\Phi_1 \ast \Phi_2: g \longmapsto \int_{G(\mathbb{A})} \Phi_1(gh^{-1})\Phi_2(h)dh$$
for $\Phi_1, \Phi_2 \in \mathcal{H}$ is called the \textit{Hecke algebra for }$G(\mathbb{A})$. Its elements are called \textit{Hecke operators}.
\end{definition}

The zero element of $\mathcal{H}$ is the zero function, but there is no multiplicative unit. For $K' \in \mathcal{V},$ we define $\mathcal{H}_{K'}$ to be the subalgebra of all bi-$K'$-invariant elements. These subalgebras have multiplicative units. Namely, the normalized characteristic function $\epsilon_{K'} := (\mathrm{vol}K')^{-1} \mathrm{char}_{K'}$ acts as the identity on $\mathcal{H}_{K'}$ by convolution. When $K'=K$ we call $\mathcal{H}_K$ the \textit{unramified part of} $\mathcal{H}$ and its elements are called \textit{unramified Hecke operators}. 

\begin{proposition}  $\mathcal{H} = \bigcup_{K' \in \mathcal{V}} \mathcal{H}_{K'}.$

\end{proposition}

The Hecke algebra $\mathcal{H}$ acts on $C^{0}(G(\mathbb{A}))$ by
$$\Phi(f): g \longmapsto \int_{G(\mathbb{A})} \Phi(h) f(gh) dh.$$

A function $f \in C^{0}(G(\mathbb{A}))$ is $\mathcal{H}-$\textit{finite} if the space $\mathcal{H}. f$ is finite dimensional.

\begin{definition} The \textit{space of automorphic forms} $\mathcal{A}$ (with trivial central character) is the complex vector space of all functions $f \in C^{0}(G(\mathbb{A}))$ which are smooth, $K-$finite, of moderate growth,  left $G(F)Z(\mathbb{A})-$invariant and $\mathcal{H}-$finite. Its elements are called automorphic forms. 

\end{definition}

\begin{lemma}[\cite{oliver-thesis} Lemma 1.3.2]  A function $f \in C^{0}(G(\mathbb{A}))$ is smooth and $K-$finite if and only if there is a $K' \in \mathcal{V}$ such that $f$ is right $K'-$invariant. 
\end{lemma}

For every subspace $V \subseteq \mathcal{A},$ let $V^{K'}$ be the subspace of all $f \in V$ that are right $K'-$invariant. By the previous lemma, functions in $\mathcal{A}^{K'}$ can be identified with the functions on $G(F)Z(\mathbb{A}) \setminus G(\mathbb{A})/K'$ that are of moderate growth.

A consequence of last lemma is:

\begin{proposition} $V = \bigcup_{K' \in \mathcal{V}} V^{K'}$ for every subspace $V \subseteq \mathcal{A}.$ 
$\hfill \square$
\end{proposition}


The fundamental fact for defining the main object of this paper is the following proposition. 

\begin{proposition}  \label{propfund} Fix $\Phi \in \mathcal{H}_{K'}$. For all  $[g] \in G(F) \setminus G(\mathbb{A}) / K',$ there is a unique set of pairwise distinct classes $[g_1], \ldots,[g_r] \in G(F) \setminus G(\mathbb{A}) / K'$ and numbers $m_1,  \ldots, m_r \in \C^{*}$ such that  
$$\Phi(f)(g) = \sum_{i=1}^{r} m_i f(g_i).$$
for all $f \in \mathcal{A}^{K'}.$

\begin{proof} The existence is as follows. Since $\Phi$ is $K'-$bi-invariant and compactly supported, it is a finite linear combination of characteristic functions on double cosets of the form $K'hK'$ with $h \in G(\mathbb{A}).$ So we may reduce the proof to $\Phi = \mathrm{char}_{K'hK'}.$ Again, since $K'hK'$ is compact, it equals the union of a finite number of pairwise distinct cosets $h_1K', \ldots, h_rK',$ and thus
$$\begin{array}{ccc}
\Phi(f)(g) & = & \hspace{-0.7cm}\int_{G(\mathbb{A})} \mathrm{char}_{K'hK'}(h') f(gh') dh' \\ 

 & = & \hspace{-0.1cm}\sum_{i=1}^{r} \int_{G(\mathbb{A})} \mathrm{char}_{h_iK'}(h') f(gh') dh' \\ 
 
 & = & \hspace{-2.1cm}\sum_{i=1}^{r} \mathrm{vol}(K')f(gh_i) \\ 

 & = & \hspace{-3cm}\sum_{i=1}^{r} m_i f(g_i)
\end{array}  $$
for any $g \in G(\mathbb{A})$, where for the last equality we have taken care of putting together values of $f$ in the same classes of $\in G(F) \setminus G(\mathbb{A}) / K'$ and throwing out zero terms. Uniqueness follows of the construction and the fact that $f \in \mathcal{A}^{K'}.$ 
\end{proof}
\end{proposition}

For $[g],[g_1], \ldots,[g_r] \in G(F) \setminus G(\mathbb{A}) / K',$ as in the last proposition, we denote  $\mathcal{V}_{\Phi,K'}([g]) := \{([g],[g_i],m_i)\}_{i=1, \ldots, r}.$

\begin{definition} \label{defgraphs} Using the notation of Proposition \ref{propfund}, we define the graph $\mathcal{G}_{\Phi,K'}$ of $\Phi$ relative to $K'$ whose vertices are 
$$\mathrm{Vert} \mathcal{G}_{\Phi,K'} = G(F) \setminus G(\mathbb{A}) / K'$$
and the oriented weighted edges 
$$\mathrm{Edge} \mathcal{G}_{\Phi,K'} = \bigcup_{[g] \in \mathrm{Vert} \mathcal{G}_{\Phi,K'}} \mathcal{V}_{\Phi,K'}([g]).$$
The classes $[g_i]$ are called the $\Phi-$neighbors of $[g]$ (relative to $K'$). 
\end{definition}

We make the following drawing conventions to  illustrate the graph of a Hecke operator: vertices are represented by 
labelled dots, and an edge $([g],[g'],m)$ together with its origin $[g]$ and its terminus $[g']$ is drawn as

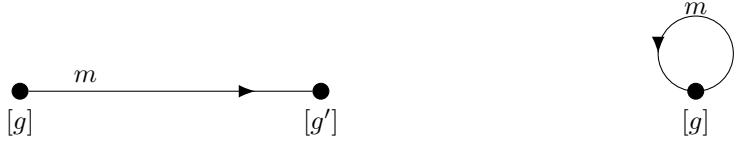
\begin{center}

\begin{figure}[h]
\centering

\begin{minipage}[b]{0.45\linewidth}
\[
  \beginpgfgraphicnamed{tikz/fig1}
  \begin{tikzpicture}[>=latex, scale=2]
        \vertex[circle,fill,label={below:$[g]$}](00) at (0,0) {};
        \vertex[circle,fill,label={below:$[g']$}](10) at (2,0) {};
    \path[-,font=\scriptsize]
    (00) edge[->-=0.8] node[pos=0.2,auto,black] {\normalsize $m$} (10)
    ;
   \end{tikzpicture}
 \endpgfgraphicnamed
\]
\end{minipage} 
\hfill
\begin{minipage}[b]{0.45\linewidth}
\[
  \beginpgfgraphicnamed{tikz/fig2}
  \begin{tikzpicture}[>=latex, scale=2]
        \vertex[circle,fill,label={below:$[g]$}](00) at (0,0) {};
  \draw[->-=0.5] (0,0.25) circle  (0.25cm) node at (0,0.55) {\normalsize $m$} ;
   \end{tikzpicture}
 \endpgfgraphicnamed
\]

\end{minipage}
\end{figure}

\end{center}

where the second figure is the case $[g]=[g'].$ 

By Proposition \ref{propfund} and the definition of the graph of $\Phi$, we have for $f \in \mathcal{A}^{K'}$ and $[g] \in G(F) \setminus G(\mathbb{A}) / K'$ that
$$\Phi(f)(g) = \displaystyle\sum_{\substack{([g],[g_i],m_i)\\ \in \mathrm{Edge}\mathcal{G}_{\Phi,K'}}} m_i f(g_i).$$
Hence one can read off the action of a Hecke operator on the value of an automorphic form from the illustration of the graph: 

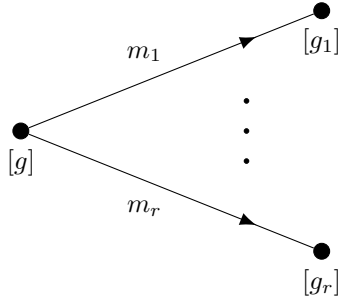
\[
  \beginpgfgraphicnamed{tikz/fig}
  \begin{tikzpicture}[>=latex, scale=2]
        \vertex[circle,fill,label={below:$[g]$}](00) at (0,0) {};
        \vertex[circle,fill,label={below:$[g_1]$}](11) at (2,0.8) {};
        \vertex[circle,fill,label={below:$[g_r]$}](10) at (2,-0.8) {};
\draw (1.5,0.2) circle (0.015cm);
\fill  (1.5,0.2) circle (0.015cm);      
\draw (1.5,0) circle (0.015cm);
\fill  (1.5,0) circle (0.015cm);  
\draw (1.5,-0.2) circle (0.015cm);
\fill  (1.5,-0.2) circle (0.015cm);   
            \path[-,font=\scriptsize]
    (00) edge[->-=0.8] node[pos=0.5,auto,black,swap] {\normalsize $m_r$} (10)
    (00) edge[->-=0.8] node[pos=0.5,auto,black] {\normalsize  $m_1$} (11)
    ;
   \end{tikzpicture}
 \endpgfgraphicnamed
\]

The next proposition is an easy consequence of the algebra structure of $\mathcal{H}_{K'}.$ 

\begin{proposition} \label{propedge}For the zero element $0 \in \mathcal{H}_{K'},$ the multiplicative unit $1 \in \mathcal{H}_{K'}$ and arbitrary $\Phi_1, \Phi_2 \in \mathcal{H}_{K'}, r \in \C^{*},$ we obtain that
\vspace{0.2cm}

\begin{enumerate}
\item[$(i)$]$\mathrm{Edge}\; \mathcal{G}_{0,K'} = \emptyset;$ \\

\item[$(ii)$]  $\mathrm{Edge}\; \mathcal{G}_{1,K'} = \big\{([g], [g],1)\big\}_{[g] \in \mathrm{Vert} \; \mathcal{G}_{1,K'} };$\\

\item[$(iii)$] $\mathrm{Edge} \; \mathcal{G}_{\Phi_1 + \Phi_2,K'} = \Big\{([g],[g'],m) \; \Big| \; { \small m = \displaystyle\sum_{\substack{([g],[g'],m')\\ \in \mathrm{Edge} \; \mathcal{G}_{\Phi_1,K'} }} m'  + \sum_{\substack{([g],[g'],m'') \\ \in \mathrm{Edge} \; \mathcal{G}_{\Phi_2,K'} }} m'' \neq 0 }\Big\};$\\

\item[$(iv)$] $\mathrm{Edge} \; \mathcal{G}_{r\Phi_1,K'} = \Big\{([g],[g'],rm) \; \big| \; ([g],[g'],m) \in \mathrm{Edge}\; \mathcal{G}_{\Phi_1,K'}\Big\},$ and\\

\item[$(v)$] $\mathrm{Edge} \; \mathcal{G}_{\Phi_1 \ast \Phi_2,K'} = \Big\{([g],[g'],m)  \; \Big|\;{\small m = \displaystyle\sum_{\substack{([g],[g''],m') \in \mathrm{Edge}\; \mathcal{G}_{\Phi_1,K'} \\ ([g''],[g'],m'') \in \mathrm{Edge}\; \mathcal{G}_{\Phi_2,K'}}} m' \cdot m'' \neq 0  }\Big\}.$
\end{enumerate}

$ \hfill \square$

\end{proposition}


\section{Graphs of unramified Hecke operators} \label{section2}

In this section, we investigate the structure around a vertex in the graph of unramified Hecke operators. The main result of this section is that the neighbors of a vertex, counted with multiplicities, correspond to the points of an appropriate Grassmannian.

\textbf{Notation:} We adopt the convention that the empty entry in a matrix means a zero entry.


In the following given an explicit description of the unramified Hecke algebra. Fix $n \geq 1$ an integer. For $x$ a place of $F,$ let $\Phi_{x,r}$ be the characteristic function of
$$K\left( \begin{array}{cc}
\pi_x I_r &  \\ 
 & I_{n-r}
\end{array}    \right)K$$  
where $I_k$ is the $k \times k$ identity matrix. Observe that $\Phi_{x,n}$ is invertible and its inverse is given by the characteristic function of $K (\pi_x I_n)^{-1} K.$ The well-known theorem about the polynomial structure of unramified Hecke algebra due by Tamagawa and Satake is:

\begin{theorem}[\cite{dennis} Chapter 12, 1.6] \label{satake}  Identifying $\epsilon_K$ with $1 \in \C$ yields
$$\mathcal{H}_K \cong \C[\Phi_{x,1}, \ldots, \Phi_{x,n}, \Phi_{x,n}^{-1}]_{x \in |X|}.$$
In particular, $\mathcal{H}_K $ is commutative. 

\end{theorem}

By Proposition \ref{propedge}, it is enough to determine the graphs for the algebra generators $\Phi_{x,1}, \ldots, \Phi_{x,n}, \Phi_{x,n}^{-1}$ with $x \in |X|$ (Theorem \ref{satake}) in order to understand the graph of any Hecke operator. We use the shorthand notation $\mathcal{G}_{x,r}$ for the graph $\mathcal{G}_{\Phi_{x,r},K}$ and $\mathcal{V}_{x,r}([g])$ for $\Phi_{x,r}-$neighborhood $\mathcal{V}_{\Phi_{x,r}, K}([g])$ of $[g]$, for $x \in |X|$ and $r=1, \ldots, n.$

The following considerations will be used in the proof of Theorem \ref{maintheorem}. For the standard Borel subgroup $B < G$ of upper triangular matrices, we have the local and global form of \textit{Iwasawa decomposition}, respectively:
$$G(F_x) = B(F_x)K_x \hspace{0.5cm}\text{ and } \hspace{0.5cm} G(\mathbb{A}) = B(\mathbb{A})K$$
with $K_x = G(\mathcal{O}_x),$ see  Proposition 4.5.2 in \cite{Bump}. 

\vspace*{0.3cm}
\noindent
\textbf{Schubert Cell Decomposition.} Let $\mathrm{Gr}(k,n)$ denote the Grassmannian that parametrizes the $k-$dimensional linear subspaces of a fixed $n-$dimensional vector space. As a set, the Grassmannian has a decomposition as the disjoint union 
$$\mathrm{Gr}(k,n) = \coprod_{\lambda \in J(k,n)} C_\lambda$$
where $J(k,n) = \{(j_1, \ldots,j_k) | 1 \leq j_1 < \cdots < j_k \leq n\},$  and
 $C_\lambda$ is the set of $n \times n-$matrices $(a_{ij})_{n \times n}$ of the following form:
\begin{itemize}
\item $a_{jj} =1$ if $j \in \lambda,$

\item $a_{ij}=0$ if $j \not\in \lambda,$ or $j<i,$ or $i \in \lambda$ and $j \not\in \lambda,$ 
\end{itemize} 
where we denote $j \in \lambda= (j_1, \ldots, j_n)$ if $j \in \{j_1, \ldots, j_n\}.$ 
The sets $C_\lambda$ are called \textit{Schubert cells} and the disjoint union above is the \textit{Schubert cell decomposition} of $\mathrm{Gr}(k,n).$ See  Theorem $11$ in \cite{hoffman-kunze} and \cite{Fulton} for more details. 




We associate with each $w \in \mathrm{Gr}(k,n)(\mathbb{F}_{q})$ a matrix $\xi_w \in G(\mathbb{A}).$ By the Schubert cell decomposition, each $w \in \mathrm{Gr}(k,n)(\mathbb{F}_{q})$ lies in precisely one Schubert cell $C_\lambda.$ Thus we can identify $w$ with a matrix $(a_{ij})_{n \times n } $ in $C_\lambda.$ We define $\xi_w = (b_{ij}) \in G(\mathbb{A})$ by just replacing $0$'s on the diagonal of $(a_{ij})_{n \times n } $ by $\pi_x.$ 




We identify $\pi_x$ with the idele $a$ where $a_y = 1$ if $y \neq x$ and $a_x= \pi_x$ in $y=x$, and we can consider $\alpha \in \kappa(x)$ as the adele whose component at $x$ is $\alpha$ and whose others components are $0.$ 

With these identifications, we can consider $\xi_w$ as an adelic matrix in $G(\mathbb{A}).$ The Schubert cell decomposition can be reformulated as follows. 

\begin{lemma} \label{lemma1}There is a bijection of  $ \mathrm{Gr}(k,n)(\kappa(x))$ with the set
$$\left\lbrace \left( \begin{small}
\begin{array}{cccc}
\epsilon_1 & b_{12} & \cdots & b_{1n} \\ 
 & \ddots &  &  \vdots\\ 
 &  &  \epsilon_{n-1}& b_{n-1 n} \\ 
 &  &  & \epsilon_n
\end{array}   \end{small} \right) \left|  
\begin{array}{c}
\epsilon_i \in \{1, \pi_x \}, \\ 
\#\{i | \epsilon_i = 1\} =k \text{ and } b_{ij} \in \kappa(x) \\ 
  \text{ with } b_{ij} = 0 \text{ if  either}\\ 
  \epsilon_j =\pi_{x}^{} \text{ or }  \epsilon_i =1
\end{array}  \right. \right\rbrace 
$$

\end{lemma}


\begin{example} For $n=2$ and $k=1$, we have
$$\xi_{w_1} = \left( \begin{array}{cc}
1 & 0 \\ 
0 & \pi_x
\end{array} \right)\hspace{0.2cm} \text{ where }\hspace{0.2cm}  w_1 = [1:0] \in \mathbb{P}^1(\kappa(x))$$
and
$$\xi_{w_2} = \left( \begin{array}{cc}
\pi_x & \ast \\ 
0 & 1
\end{array} \right)\hspace{0.2cm}  \text{ where } \hspace{0.2cm} w_2 = [\ast : 1] \in \mathbb{P}^1(\kappa(x)), $$
which are the same matrices as considered by Lorscheid in \cite{oliver-thesis}. \end{example}




As the next step, we describe a method to distinguish the double cosets $K_x \setminus G(F_x) / K_x$, which will be useful in the proof of the next lemma. Let
$$\begin{array}{cccc} \bigwedge^k:  &
\mathrm{GL}_n (F_x) & \longrightarrow & \mathrm{GL}_{{n}\choose{k}}(F_x) \\[0.3cm] 
& g & \longmapsto & \wedge^k g
\end{array} $$
be the $k-$th exterior power representation, i.e.\ the entries of $\wedge^k g$ are the $k \times k-$minors of $g.$ Let $I(\wedge^k g)$ be the fractional ideal of $\mathcal{O}_x$ generated by the entries in $\wedge^k g$, i.e.\ the fractional ideal generated by these minors. It is clear that $\wedge^k K_x$ is a subset of $\mathrm{GL}_{{n}\choose{k}}(\mathcal{O}_x),$ thus $I(\wedge^k g)$ is invariant under left and right multiplication by $K_x.$ Also note that $I(\wedge^1 g) = I(g)$ is the fractional ideal generated by the entries of $g$ and that $I(\wedge^n g)$ is generated by  determinant of $g$. Thus the following lemma.

\begin{lemma} \label{lemma2} Let $g_1, g_2 \in \mathrm{GL}_{n}(F_x)$. If the double cosets $ K_x g_1 K_x$ and $ K_x g_2 K_x$ are equal, then $I(\wedge^k g_1) = I(\wedge^k g_2)$ for $k=1, \ldots,n. \hfill \square$ 

\end{lemma}

The next lemma is the key observation for the proof of the main theorem of this section. 

\begin{lemma} \label{hellslemma} There is a decomposition of sets 
$$K \left( \begin{array}{cc}
\pi_x I_r &  \\ 
 & I_{n-r}
\end{array} \right) K = \coprod_{w \in \mathrm{Gr}(n-r,n)(\kappa(x))} \xi_w K.$$

\begin{proof} Consider the map 
$$\begin{array}{ccc}
K \times K & \longrightarrow   & G(\mathbb{A}) \\ 
(k_1,k_2) & \longmapsto & k_1 \left( \begin{array}{cc}
\pi_x I_r &  \\ 
 & I_{n-r}
\end{array} \right) k_2
\end{array} $$
which is continuous because it is induced by the group multiplication. By Tychonoff's theorem, $K \times K$ is compact. As the continuous image of the compact set $K \times K$, the double coset
$$K \left( \begin{array}{cc}
\pi_x I_r &  \\ 
 & I_{n-r}
\end{array} \right) K$$
is compact. Since $K$ is open, the quotient
$$K \left( \begin{array}{cc}
\pi_x I_r &  \\ 
 & I_{n-r}
\end{array} \right) K \Bigr/ K $$
is finite. Therefore,
$$K \left( \begin{array}{cc}
\pi_x I_r &  \\ 
 & I_{n-r}
\end{array} \right) K = \bigsqcup_{\xi \in I} \xi K$$
where $I$ is a finite set of coset representatives. 

Thus we have to show that $I = \big\{\xi_w | \  \ w \in \mathrm{Gr}(n-k,n)(\kappa(x))\big\}$ as in the Lemma  \ref{lemma1}. By the Iwasawa decomposition, we can take $\xi \in I$ to be upper triangular. The question can be solved component-wise at each place $y \in |X|.$

If $y \neq x,$ then $K_y I_n K_y = K_y.$ Thus we can use
$$w= \left( \begin{array}{c}
 I_{n-r}  \\ 
\\ \end{array}\right)_{n \times (n-r)}$$
and obtain
$$\xi_w = \left( \begin{array}{cc}
I_{n-r} &  \\ 
	 & \pi_x I_{r}
\end{array}  \right) \hspace{0.2cm} \text{ with } \hspace{0.2cm} (\xi_w)_y = I_n,$$
as desired. 

If $y=x,$ observe that 
$$K_x \left( \begin{array}{cc}
\pi_x I_r &  \\ 
 & I_{n-r}
\end{array} \right) K_x \subseteq \mathrm{Mat}(\mathcal{O}_x).$$
Thus $\xi_x$ is an upper triangular matrix with entries in $\mathcal{O}_x.$ Let
$$g = \left( \begin{small}\begin{array}{cc}
\pi_x I_{r} &  \\ 
 & I_{n-r}
\end{array} \end{small}\right).$$
Lemma \ref{lemma2} states that $I(\wedge^k g) = I(\wedge^k \xi_x)$ for every $k=1, \ldots,n.$ By a simple calculation, we have $I(g)=  \mathcal{O}_x, \   \ I(\wedge^2 g) = \mathcal{O}_x, \   \ \ldots \   \ , \   \  I(\wedge^{n-r-1}g)= \mathcal{O}_x, \   \ I(\wedge^{n-r} g)= \mathcal{O}_x,  \    \ I(\wedge^{n-r+1} g) = (\pi_x), \   \ \ldots \   \ , I(\wedge^n g) = (\pi_{x}^{r}).$
Let
$$\xi_x = \left( \begin{small}\begin{array}{ccc}
\epsilon_1 &  & b_{ij} \\ 
 & \ddots &  \\ 
 & 	 & \epsilon_n
\end{array} \end{small}\right) $$
where $\epsilon_i = u_i \pi_{x}^{l_i}, l_i \geq 0$ and $u_i \in \mathcal{O}_{x}^{*}.$ By right multiplication with
$$\left(\begin{small} \begin{array}{ccc}
u_{1}^{-1} &  &  \\ 
 & \ddots &  \\ 
 & 	 & u_{n}^{-1}
\end{array} \end{small}\right) \in K,$$
we can assume that $\epsilon_i = \pi_{x}^{l_i}, l_i \geq 0,$ thus
$$\xi_x = \left( \begin{small}\begin{array}{ccc}
\pi_{x}^{l_1} &  & b_{ij} \\ 
 & \ddots &  \\ 
 & 	 & \pi_{x}^{l_n}
\end{array}\end{small} \right) \text{ with } b_{ij} \in \mathcal{O}_x.$$
As $I(\wedge^n \xi_x) = I(\wedge^n g) = (\pi_{x}^{r}),$ we have $l_1 + \cdots + l_n = r.$ 

Suppose that there is some $l_j > 1.$ As $l_i \geq 0,$  there are at least $n-r+1$ $l_i's$ equal to zero. Thus $I(\wedge^{n-r+1} \xi_x) = \mathcal{O}_x,$ but $I(\wedge^{n-r+1} g) = (\pi_{x}^{})$,  a contradiction. Therefore, all $l_i \in \{0,1\}$ and  $\#\{i \in \{1, \ldots,n\}\big| l_i = 1\} =r.$ This show that the diagonal of $\xi_x$ is of the predicted shape. 

We are left with showing that the $b_{ij}$ are as claimed. First we can assume that $b_{ij} =0$ if $\epsilon_i =1.$ This follows since we can change the representative of $\xi K$ by means of the product
$$\left(\begin{small} \begin{array}{cccccc}
\epsilon_1 & b_{12} &  & \cdots &  & b_{1n}\\ 
 & \ddots &  &  &  & \vdots \\ 
 &  & 1 & b_{ii+1} & \cdots & b_{in}\\ 
 &  &  & \ddots &  &\vdots \\ 
 &  &  &  & \ddots & \vdots \\
 & & & & & \epsilon_n
\end{array} \end{small} \right) \cdot
\left( \begin{small}\begin{array}{cccccc}
1 &  &  &  &  &  \\ 
 & \ddots &  &  &  &  \\ 
 &  & 1 & -b_{i i+1} & \cdots & - b_{in} \\ 
 &  &  & \ddots &  &  \\ 
 &  &  &  & \ddots &  \\ 
 &  &  &  &  & 1 \end{array} \end{small} \right) .$$ 
Doing this step for every $i \in J$, we obtain a matrix $\xi_x$ which is zero in entries on right and left of $\epsilon_i = 1$. Observe that these multiplications change the original $b_{ij}  's$, for convenience we shall continue using the same notation for those elements.

If $\epsilon_i = \epsilon_j = \pi_{x}$, then $b_{ij} \in I(\wedge^{n-r+1} \xi_x).$ Indeed, let $B_{ij}$ the $(n-r+1) \times (n-r+1)$ submatrix of $\xi_x$ by removing the $r-1$ columns $\neq j$ with $\pi_x$ in the diagonal and removing the $r-1$ rows $\neq i$ with $\pi_x$ in the diagonal. Observe that $b_{ij}$ is an entry of $B_{ij}$ and  the column of $B_{ij}$ which contains $b_{ij}$ has all entries, except $b_{ij},$ equal to zero. Moreover, when we remove the row and the column of $B_{ij}$ which contains $b_{ij},$ we get the identity matrix. Thus, from Laplace's formula, $b_{ij} = \pm \det(B_{ij})$ and hence $b_{ij} \in I(\wedge^{n-r+1} \xi_x).$

Since $(\pi_x) = I(\wedge^{n-r+1} g) = I(\wedge^{n-r+1} \xi_x)$, we have $\pi_{x}^{} \mid b_{ij}$ if $\epsilon_i = \epsilon_j = \pi_{x}^{}.$ Therefore we can eliminate these entries by the following multiplication
$$ \left(\begin{small} \begin{array}{ccccccc}
\epsilon_1 &  &  &  &  &  &  \\ 
 & \ddots &  &  &  &  &  \\ 
 & & \pi_{x}^{}&  & \pi_{x}^{} \tilde{b_{ij}} &  &  \\ 
 &  &  & \ddots &  &  &  \\ 
 &  &  &  & \pi_{x}^{} &  &  \\ 
 &  &  &  &  & \ddots&  \\ 
 &  &  &  &  &  & \epsilon_n
\end{array}\end{small} \right) \cdot 
\left(  \begin{small}\begin{array}{ccccccc}
1 &  &  &  &  &  &  \\ 
 & \ddots &  &  &  &  &  \\ 
 & & 1&  &  -\tilde{b_{ij}} &  &  \\ 
 &  &  & \ddots &  &  &  \\ 
 &  &  &  & 1 &  &  \\ 
 &  &  &  &  & \ddots&  \\ 
 &  &  &  &  &  & 1
\end{array} \end{small}\right) $$ 
where $\pi_{x}^{} \tilde{b_{ij}} = b_{ij}.$

Finally, for $b_{ij} = \sum_{\ell \geq 0} b_{ij}(\ell) \pi_{x}^{\ell}$ where $\epsilon_j=1$ and $\epsilon_i = \pi_{x}^{}$, we multiply $\xi_x$ with the  matrix $\lambda = (c_{rs})$ from the right, where $c_{rr}=1$, and for $r \neq s$, $c_{rs}=0$ unless $r=i,s=j,$ in which case $c_{ij} =  -\sum_{\ell \geq 1} b_{ij}(\ell) \pi_{x}^{\ell -1}.$ This allow to consider $b_{ij}$ equal to $b_{ij}(0) \in \kappa(x)$, therefore $\xi_x$ is as claimed  in the lemma \ref{lemma1}.
\end{proof}
\end{lemma}

\begin{theorem} \label{maintheorem} The $\Phi_{x,r}-$neighbors of $[g]$ are the classes $[g \xi_w]$ with $\xi_w$ as in the previous lemma and the multiplicity of an edge of  $[g]$ to $[g']$ equals the number of $w \in \mathrm{Gr}(n-r,n)(\kappa(x))$ such that $[g \xi_w] = [g'].$ The multiplicities of the edges originating in $[g]$ sum up to $\#\mathrm{Gr}(n-r,n)(\kappa(x)).$

\begin{proof} Let $\Phi := \mathrm{char}(KhK)$, for $h \in G(\mathbb{A}).$ Since $KhK$ is compact and $K$ is open, $Kh K$ is equals to the union of a finite number of pairwise distinct cosets $h_1K, \ldots, h_rK.$ 
Hence, by Proposition \ref{propfund} and definition of $\mathcal{G}_{\Phi,K}$ the neighbors of $[g]$ are $[gh_i], i=1, \ldots,r.$ Therefore the theorem follows from preceding lemma.   
\end{proof}
\end{theorem}

 



\begin{corollary} \label{corollaryofmaintheorem}The edges set of the graph of $\Phi_{x,n}$ is  
$$\mathrm{Edge} \mathcal{G}_{x,n} = \big\{([g], [g (\pi_x I_n)],1)\big\}_{[g] \in G(F)\setminus G(\mathbb{A}) / K }. $$

\begin{proof} Observe that $K (\pi_{x} I_n) K = (\pi_{x} I_n) K.$
\end{proof}

\end{corollary}

In the calculation of the graph of $\Phi_{x,r},$ the following well-known formula is very useful. 





\begin{lemma} $\# \mathrm{Gr}(k,n)(\mathbb{F}_q) = \dfrac{(q^n -1)(q^{n-1} - 1)(q^{n-2} - 1) \cdots (q^{n-k+1} - 1)}{(q^k -1)(q^{k-1} - 1)(q^{k-2}-1) \cdots (q -1)}.$

\end{lemma}

\section{Geometry of graphs of unramified Hecke operators} \label{section3}

Recall that $X$ is the geometrically irreducible smooth projective curve over $\Fq$ whose function field is $F.$ A well-known theorem by Weil states that $G(F)\setminus G(\mathbb{A}) /K$ stays in bijection to the set $\mathrm{Bun}_n X$ of isomorphism classes of rank-$n$ vector bundles on $X.$ This allows us to give a interpretation of $\mathcal{G}_{x,r}$ in geometric terms. We begin with a review of Weil's theorem.


The bijection
$$\begin{array}{ccccc}
F^{*}\setminus \mathbb{A}^{*} / \mathcal{O}_{\mathbb{A}}^{*} & = \mathrm{Cl}F & \stackrel{1:1}{\longleftrightarrow} & \mathrm{Pic}X = & \mathrm{Bun}_{1}X \\
{[a]} & & \longmapsto & & \Line_{a}  \end{array} $$ 
where $\Line_{a} = \Line_{D}$ if $D$ is the divisor determined by $a$, generalises to all vector bundles as follows, cf.  Lemma $3.1$ in \cite{frenkel2} and  $2.1$ in \cite{dennis}. 
A rank $n$ bundle $\E$ can be described by choosing bases 
$$\E_{\eta} \cong \mathcal{O}_{X, \eta}^{n} = F^{n} \quad \text{ and } \quad \E_x \cong \mathcal{O}_{X,x}^{n} = (\mathcal{O}_x \cap F)^{n}$$
for all stalks, where $\eta$ is the generic point of $X,$ and the inclusion maps $\E_x \hookrightarrow \E_{\eta}$ for all $x \in |X|.$ After tensoring with $F_x$, for each $x \in |X|$ we obtain
$$F_{x}^{n} \cong \mathcal{O}_{x}^{n} \otimes_{\mathcal{O}_{x}} F_{x} \cong \E_{x} \otimes_{\mathcal{O}_{X,x}} \mathcal{O}_{x} \otimes_{\mathcal{O}_{x}} F_{x} \cong \E_{x} \otimes_{\mathcal{O}_{X,x}} F \otimes_{F} F_{x} \cong F^{n} \otimes_{F} F_{x} \cong F_{x}^{n} $$
which yields an element $g$ of $G(\mathbb{A}).$ 

A change of bases for $\E_{\eta}$ and $\E_x$ corresponds to multiplying $g$ by an element of $G(F)$ from the left and by an element of $K$ from the right, respectively.

Since the inclusion $F \subset F_x$ is dense for every place $x$, and $G(\mathcal{O}_{\mathbb{A}})$  is open in $G(\mathbb{A})$,
every class in $G(F)\setminus  G(\mathbb{A})/G(\mathcal{O}_{\mathbb{A}})$ is represented by a $g= (g_x) \in G(\mathbb{A})$ such that
$g_x \in G(F)$ for all places $x$. This means that the above construction can be reversed. Weil's theorem asserts the following.


\begin{theorem}[\cite{frenkel2} Lemma 3.1]  \label{weil} For every $n \geq 1,$ the above construction yields a bijection
$$\begin{array}{ccc}
G(F) \setminus G(\mathbb{A}) / K & \stackrel{1:1}{\longleftrightarrow} & \mathrm{Bun}_n X. \\ 
g & \longmapsto & \E_g
\end{array}   $$

\end{theorem}

The last theorem identifies the set of vertices of $\mathcal{G}_{x,r}$ with  geometric objects $\mathrm{Bun}_n X.$ The next task is to describe the edges of $\mathcal{G}_{x,r}$ in geometric terms. We say that two exact sequences of sheaves
$$0 \longrightarrow \F_1 \longrightarrow \F \longrightarrow \F_{2} \longrightarrow 0 \hspace{0.2cm}\text{ and } \hspace{0.2cm} 0 \longrightarrow \F_{1}' \longrightarrow \F \longrightarrow \F_2' \longrightarrow 0$$  
are \textit{isomorphic with fixed $\F$} if there are isomorphism $\F_1 \rightarrow \F_1'$ and $\F_2 \rightarrow \F_2'$ such that 
$$\xymatrix@R5pt@C7pt{ 0 \ar[rr] && \F_1 \ar[rr] \ar[dd]^{\cong} && \F \ar[rr] \ar@{=}[dd] && \F_2 \ar[rr] \ar[dd]^{\cong} && 0 \\
&& && && && \\
0 \ar[rr] && \F_1' \ar[rr] && \F \ar[rr] && \F_2' \ar[rr] && 0 } $$
commutes. 
Let $\mathcal{K}_{x}$ be the torsion sheaf that is supported at $x$ and has stalk $\kappa(x)$ at $x,$ i.e. the skyscraper torsion sheaf at $x$. Fix $\E \in \mathrm{Bun}_n X$. For $r \in \{1, \ldots,n\},$ and $\E' \in \mathrm{Bun}_n X$ we define $m_{x,r}(\E,\E')$ as the number of isomorphism classes of exact sequences
$$0 \longrightarrow \E'' \longrightarrow \E \longrightarrow \mathcal{K}_{x}^{\oplus r} \longrightarrow 0$$
with fixed $\E$ and with $\E'' \cong \E'.$

\begin{definition} \label{defconection} Let $x \in |X|$. For a vector bundle $\E \in \mathrm{Bun}_n X$ we define
$$\mathcal{V}_{x,r}(\E) := \{(\E, \E', m) | m = m_{x,r}(\E,\E') \neq 0\},$$
and we call $\E'$ a $\Phi_{x,r}$-neighbor of $\E$ if $m_{x,r}(\E,\E') \neq 0$, and $m_{x,r}(\E,\E')$ its multiplicity.
\end{definition}

We will show that this concept of neighbors is the same as the one defined for classes in $G(F)\setminus G(\mathbb{A}) / K$ in Definition \ref{defgraphs}. According to Theorem \ref{maintheorem}, the $\Phi_{x,r}-$neighbors of a class $[g] \in G(F) \setminus G(\mathbb{A}) / K$ are of the form $[g\xi_w]$ for a $w \in \mathrm{Gr}(n-r,n)(\kappa(x)).$

\begin{lemma}\label{neighbourslemma} For every $x \in |X|,$ the map 

$$\begin{array}{ccc}
\mathcal{V}_{x,r}([g]) & \longrightarrow & \mathcal{V}_{x,r}([\E_g]) \\ 
([g],[g'],m) & \longmapsto & (\E_g,\E_{g'},m)
\end{array} $$
is a well-defined bijection. 

\begin{proof} Fix a base $(\E_g)_y \stackrel{\sim}{\rightarrow} \mathcal{O}_{X,y}^{\oplus n},$ for each $y \in |X|.$ Note that by the definition of $\xi_w$, multiplying an element of $\mathcal{O}_{X,y}^{\oplus n}$ with the component $(\xi_w)_y$ from the right yields an  element of $\mathcal{O}_{X,y}^{\oplus n}.$ Thus we obtain an exact sequence of $\mathbb{F}_q-$modules
$$0 \longrightarrow \prod_{y \in |X|} \mathcal{O}_{X,y}^{\oplus n} \stackrel{\xi_w}{\longrightarrow}\prod_{y \in |X|} \mathcal{O}_{X,y}^{\oplus n} \longrightarrow \kappa(x)^{\oplus r} \longrightarrow 0$$
where $r$ is the number of $\pi_x$ in the diagonal of $(\xi_w)_x.$ The correspondence in  Theorem \ref{weil} implies the following exact sequence of sheaves
$$0 \longrightarrow \E_{g \xi_w} \longrightarrow \E_g \longrightarrow \mathcal{K}_{x}^{\oplus r} \longrightarrow 0.$$
This maps $w \in \mathrm{Gr}(n-r,n)(\kappa(x))$ to the isomorphism class of 
$$\big(0 \rightarrow \E_{g \xi_w} \rightarrow \E_g \rightarrow \mathcal{K}_{x}^{\oplus r} \rightarrow 0\big)$$ 
with fixed $\E_g.$

On the other hand, an isomorphism class of exact sequences
\begin{equation}\label{eq3.1}
\big(0 \longrightarrow \E' \longrightarrow \E \longrightarrow \mathcal{K}_{x}^{\oplus r} \longrightarrow 0\big)  
\end{equation}
of sheaves yields an isomorphism class of exact sequence of $\mathcal{O}_{X,x}$-modules 
$$0 \longrightarrow \E_x' \longrightarrow \E_x \longrightarrow \kappa(x)^{\oplus r} \longrightarrow 0$$
by taking stalks.
After tensoring by $\mathcal{O}_{X,x}\big/m_{X,x}$, where $m_{X,x}$ is the maximal ideal of $\mathcal{O}_{X,x},$ we obtain an exact sequence of $\kappa(x)-$vector spaces
$$ \E_x' \otimes \frac{\mathcal{O}_{X,x}}{m_{X,x}} \longrightarrow \E_x \otimes \frac{\mathcal{O}_{X,x}}{m_{X,x}} \longrightarrow \kappa(x)^{\oplus r} \longrightarrow 0.$$
Since 
$$\E_x' \otimes \frac{\mathcal{O}_{X,x}}{m_{X,x}} \cong \frac{\E_x'}{m_{X,x} \E_x'} \hspace{0.5cm} \text{ and } \hspace{0.5cm} \E_x \otimes \frac{\mathcal{O}_{X,x}}{m_{X,x}} \cong \frac{\E_x}{m_{X,x} \E_x}$$
(see \cite{atiyah} page 31) we have
$$ \frac{\E_x'}{m_{X,x} \E_x'} \longrightarrow \frac{\E_x}{m_{X,x} \E_x} \longrightarrow \kappa(x)^{\oplus r} \longrightarrow 0.$$
Let $w_x$ be the image of $\E_x'\big/m_{X,x} \E_x'$ in the above sequence. Then 
$$0 \longrightarrow w_x \longrightarrow \frac{\E_x}{m_{X,x} \E_x} \longrightarrow \kappa(x)^{\oplus r} \longrightarrow 0$$
and $\dim(w_x) = n-r.$ Hence we have associated with the exact sequence  \eqref{eq3.1} an element $w_x$ in $\mathrm{Gr}(n-r,n)(\kappa(x)).$ 

As we have chosen a bases for the stalk at $x$, and the identity map for all others closed points, these maps are inverse one of each other. Therefore, follows the desired bijection.
\end{proof}

\end{lemma}

Theorem \ref{weil} and Lemma \ref{neighbourslemma} imply: 

\begin{theorem} \label{theoremcorrespondence} Let $x \in |X|.$ The graph $\mathcal{G}_{x,r}$ of $\Phi_{x,r}$ is described in geometric terms as 
$$\mathrm{Vert}\; \mathcal{G}_{x,r} = \mathrm{Bun}_n X \hspace{0.3cm}\text{ and } \hspace{0.3cm} \mathrm{Edge} \;\mathcal{G}_{x,r} = \coprod_{\E \in \mathrm{Bun}_n X} \mathcal{V}_{x,r}(\E).$$ 
 \hspace{12.3cm} $\square$
\end{theorem}

\begin{proposition}\label{propend} Let $\Line_1, \ldots, \Line_n$ be  invertible sheaves on $X$. Then for every $S \subseteq \{1, \ldots, n\},$ with $\#S=r$
$$m=m_{x,r}(\E,\E') \neq 0, \text{ i.e. } (\E,\E',m) \in \mathcal{V}_{x,r}(\E)$$
where $ \E =  \Line_1 \oplus \cdots \oplus \Line_n$ and $\E' = \bigoplus_{j\in S}\Line_j( - x) \oplus \bigoplus_{j \not\in S} \Line_j.   $ 

\begin{proof}Consider the exact sequence 
$$0 \longrightarrow \Line(-x) \longrightarrow \mathcal{O}_X \longrightarrow \mathcal{K}_x \longrightarrow 0.$$
Tensoring the above sequence with the locally free sheaf $\Line_j$, does neither affect the sheaf $\mathcal{K}_x$ nor the exactness. Thus we get 
$$0 \longrightarrow \Line_j(-x) \longrightarrow \Line_j \longrightarrow \mathcal{K}_x \longrightarrow 0.$$
Combining these exact sequences for $j \in S$ with
$$0 \longrightarrow \Line_i \longrightarrow \Line_i \longrightarrow 0 \longrightarrow 0$$
for  $i \not\in S$ shows that
$$0 \longrightarrow \bigoplus_{j\in S}\Line_j( - x)  \oplus \bigoplus_{j \not\in S} \Line_j \longrightarrow \Line_1 \oplus \cdots \oplus \Line_n \longrightarrow \mathcal{K}_{x}^{\oplus r} \longrightarrow 0$$
is an exact sequence. 
\end{proof}

\end{proposition}



\section{The $\delta-$invariant }  \label{section4}

Let $\E$ be a locally free sheaf and $\E'$ a subsheaf. Note that the quotient $\E/\E'$ is not necessarily locally free. We will call $\E'$ a \textit{subbundle} if the quotient $\E/\E'$ is still a vector bundle i.e.\ a locally free sheaf. 

Let $\overline{\E'} := \ker\big(\E \rightarrow \E/\E' \rightarrow (\E/\E')\big/ (\E/\E')_{\mathrm{tors}}\big)$, where $(\E/\E')_{\mathrm{tors}}$ denotes the torsion subsheaf of $\E/\E'$, i.e.\ $(\E/\E')\big/(\E/\E')_{\mathrm{tors}}$ is torsion free. By definition $\overline{\E'}$ is a subbundle of $\E$, and $\E' \subset \overline{\E'}.$ It is easy to see that $\overline{\E'}/\E' \cong (\E/\E')_{\mathrm{tors}}.$ We recall that, given $\E \in \mathrm{Bun}_n X$, $\deg(\E):= \deg(\wedge^n \E)$, where $\wedge^n \E \in \mathrm{Pic}X.$ Hence  
$$\mathrm{rk}(\overline{\E'}) = \mathrm{rk}(\E') \text{ and } \deg(\overline{\E'}) \geq \deg(\E').$$
Therefore we can extend a subsheaf $\E'$ of $\E$ to a subbundle $\overline{\E'}$ of $\E.$

In what follows, $k$-(sub)bundle means a (sub)bundle of rank $k$. For a subbundle $\mathcal{E'}$ of a $n$-bundle $\E$, we define 
$$\delta(\mathcal{E}', \mathcal{E}) := \mathrm{rk}(\E) \deg(\mathcal{E}') - \mathrm{rk}(\E') \deg(\mathcal{E})$$
and for $k=1, \ldots, n-1,$
$$\delta_k(\mathcal{E}) := \sup_{\substack{\mathcal{E'} \hookrightarrow \mathcal{E} \\ k-subbundle} } \delta(\mathcal{E}', \mathcal{E}).$$
Moreover, we define 
$$\delta(\E) := \max\{\delta_1(\E), \ldots, \delta_{n-1}(\E)\}.$$


\begin{definition} We say that $\E' \hookrightarrow \E$ is a \textit{maximal $k-$subbundle} if $\delta(\E', \E) = \delta_{k}(\E)$ and that $\E'$ is a \textit{maximal subbundle} if $\delta_{}(\E',\E)=\delta(\E).$
\end{definition}

\begin{remark} A vector bundle $\E$ is (semi)stable if and only if 
$$\delta(\E)<0 \   \quad \ ( \delta(\E) \leq 0).$$

\end{remark}

\begin{lemma}\label{lemma4.3} Let $\mathcal{E}$ be a locally free
sheaf of rank $n$ on $X$. Suppose that $H^0(X, \mathcal{E})$ is non-trivial. Then there is an exact sequence
$$0 \longrightarrow \Line \longrightarrow \mathcal{E} \longrightarrow \mathcal{E}'' \longrightarrow 0$$
with $\mathcal{E}''$ locally free of rank $r-1$, $\Line$ locally free of rank $1$ and $\deg(\Line) \geq 0$.

\begin{proof} Let $s \in H^0(X,\E)$ be non-zero. The global section $s$ defines a map $\mathcal{O}_X \rightarrow \E$   
by 
$$ \begin{array}{ccc}
\mathcal{O}_X(U) & \longrightarrow & \E(U), \\ 
1 & \longmapsto & s|_{U}
\end{array} $$
This map is injective and yields the short exact sequence 
$$0 \longrightarrow \mathcal{O}_X \longrightarrow \E \longrightarrow \E/s\mathcal{O}_X \longrightarrow 0$$
Since we do not know if $\F= \E/s\mathcal{O}_X$ is locally free,  take $\E'' = \F/\F_{tors}$, which is a locally free sheaf of rank $r-1$ where $\F_{tors}$ is the subsheaf of torsion elements. Consider the surjective map
$$\psi : \E \longrightarrow \E'' $$
and define $\Line = \ker(\psi).$ By construction, we have a short exact sequence
$$0 \longrightarrow \Line \longrightarrow \E \longrightarrow \E'' \longrightarrow 0$$
with $\Line$ locally free of rank $1$ and $\mathcal{E}''$ locally free of rank $r-1$.

Concerning the degree of $\Line$, note that 
$$\deg(\Line) = \deg(\E) - \deg(\E'') = \deg(\E) - \deg(\F) + \deg(\F_{tors}) = \deg(\F_{tors}) \geq 0 $$
where the last inequality holds for every torsion sheaf. 
\end{proof}

\end{lemma}

\begin{proposition} For every rank-$n$ bundle $\mathcal{E},$ 
$$-ng \leq \delta(\mathcal{E}) < \infty$$

\begin{proof} Recall that for coherence sheaves $\mathcal{F}, \mathcal{G} \in \mathrm{Coh}(X)$, 
$$\deg(\mathcal{F} \otimes \mathcal{G}) = \mathrm{rk}(\mathcal{F})\deg(\mathcal{G}) + \deg(\mathcal{F})\mathrm{rk}(\mathcal{G}) \quad \text{ and } \quad \mathrm{rk}(\mathcal{F} \otimes \mathcal{G}) = \mathrm{rk}(\mathcal{F}) \mathrm{rk}(\mathcal{G}).$$
Let $\mathcal{L}$ be an invertible sheaf and $\E' \subseteq \mathcal{E}$ a subbundle. Then
$$\begin{array}{ccc}
\delta(\Line \otimes \E', \Line \otimes \E) & = & n\big( \deg(\E')+ \mathrm{rk}(\E') \deg(\Line)\big) - \mathrm{rk}(\E')\big(\deg(\E)  + \deg(\Line) n\big) \\ 
 & = & \hspace{-4.8cm} n \deg(\E')-\mathrm{rk}(\E') \deg(\E) \\ 
 & = &\hspace{-7.5cm}\delta(\E', \E). \\ 
\end{array} $$
Thus $\delta(-,-)$ is invariant under tensoring with line bundles. As $\deg(\Line \otimes \E) = \deg(\E) + n \deg(\Line)$, the multiplication by a line bundle changes the degree of $\E$ by a multiple of $n.$ Thus we can assume that
$$n(g-1) < \deg(\E) \leq ng.$$ 
By the Riemann-Roch theorem, 
$$\dim H^0(X, \E) = \deg(\E) + n(1-g) + \dim H^1(X,\E) \geq \deg(\E) + n(1-g) > 0.$$
By Lemma \ref{lemma4.3}, there is a line subbundle $\Line $ of $\E$ with $\deg(\Line) \geq 0.$ Thus 
$$\delta(\Line, \E)= \mathrm{rk}(\E)\deg(\Line) - \deg(\E) \geq -\deg(\E)  \geq  -ng, $$
which implies $ \delta(\E) \geq -ng.$ This establishes the first inequality.

Let $\E' \subseteq \E$ be a subbundle. Again by Riemann-Roch, 
$$\deg(\E')\leq \mathrm{rk}(\E')(g-1) + \dim H^0(X,\E') \leq \mathrm{rk}(\E)(g-1) + \dim H^0(X, \E).$$
This proves the second inequality.
\end{proof}
\end{proposition}

This proposition generalizes \cite[Prop. 6, p. 100]{serre}.
We derive some immediate consequences from the proof of the above proposition. 

\begin{corollary} The same is true in the last proposition for $\delta_1(\E)$, i.e., 
\[ \ -ng \ \leq \ \delta_1(\E)\ < \ \infty. \] $\hfill \square $
\end{corollary}

\begin{corollary} For every $k=1, \ldots, n-1, \   \ \delta_k(\E) < \infty$. $\hfill \square$
\end{corollary}

\begin{corollary} If $g=0,1$, then $ \delta_k(\E)$ is at most $n\dim H^0(X,\E)$, for every $k=1,\ldots, n-1$. In particular, $\delta(\E) \leq n \dim H^0(X,\E).$

\begin{proof} Let $\E'$ be a subbundle of $\E.$ As in the proof of last proposition  $\deg(\E') \leq \dim H^0(X,\E)$, and we can assume $\deg(\E) \geq 0.$ Thus 
$$\delta(\E',\E) = n \deg(\E') - \deg(\E) \mathrm{rk}(\E' ) \leq n \deg(\E') \leq n \dim H^0(X,\E).$$\end{proof}
\end{corollary}

\begin{lemma} Let $\E$ be a locally free and (semi)stable sheaf over $X$ such that 
\[ \deg(\E) \geq \mathrm{rk}(\E)(2g - 2),  \quad (\deg(\E) > \mathrm{rk}(\E)(2g - 2)), \] 
then $H^1(X,\E)=0.$

\begin{proof} We will prove by contradiction that $H^1(X,\E)=0.$ Suppose that $H^{1}(X,\E)\neq 0.$ Using Serre duality, 
$$H^{1}(X, \E) \cong H^{0}(X, \E^{\vee} \otimes \omega_X) \cong \mathrm{Ext}^{0}(\mathcal{O}_{X}, \E^{\vee}\otimes \omega_X) \cong \mathrm{Ext}^{0}(\E,\omega_X) = \mathrm{Hom}(\E,\omega_X).$$  
Hence there is a non-zero morphism $\varphi: \E \rightarrow \omega_X.$ Let $\mathcal{K}$ be the subbundle of $\E,$ which extends $\ker(\varphi)$. Since for line bundles $\Line, \Line'$ we have $\mathrm{Hom}(\Line,\Line') =0$ if $\deg(\Line) > \deg(\Line'),$ it is clear that $\deg(\mathrm{im}(\varphi)) \leq \deg(\omega_X).$ Thus  
$$\deg(\mathcal{K}) \geq \deg(\ker(\varphi)) \geq \deg(\E) - \deg(\omega_X) = \deg(\E) - (2g-2).$$
By (semi)stability of $\E$, 
$$\frac{\deg(\E) - (2g-2)}{\mathrm{rk}(\E)-1} \leq \frac{\deg(\mathcal{K})}{\mathrm{rk}(\E)-1}  \leq  \frac{\deg(\E)}{\mathrm{rk}(\E)}$$
where we have a strict inequality if $\E$ is stable. Therefore, $\deg(\E) < \mathrm{rk}(\E)(2g-2)$ in the stable case and $\deg(\E) \leq \mathrm{rk}(\E)(2g-2)$ in the semistable case. A contradiction with the degree of $\E.$ We conclude that $H^1(X,\E)=0.$
\end{proof}
\end{lemma}

\begin{proposition} Let $\E'$ be a subbundle of the bundle $\E$ on $X$, such that $(\E/\E')^{\vee} \otimes \E'$ is (semi)stable. If 
$$\delta(\E',\E) \geq \mathrm{rk}(\E')(\mathrm{rk}(\E)-\mathrm{rk}(\E'))(2g-2)$$
(with strict inequality in the stable case), then $\E \cong \E' \oplus \E/\E'.$

\begin{proof} We have an exact sequence
$$0 \longrightarrow \E' \longrightarrow \E \longrightarrow \E/\E' \longrightarrow 0.$$
Since, $\mathrm{Ext}^{1}(\E/\E', \E') \cong H^1(X, (\E/\E')^{\vee} \otimes \E')$ and
$$\begin{array}{ccc}
\deg((\E/\E')^{\vee} \otimes \E' ) & = & \mathrm{rk}(\E) ( \deg(\E') - \deg(\E)) + (\mathrm{rk}(\E)- \mathrm{rk}(\E'))(\deg(\E)) \\ 
 & = & \hspace{-3.4cm} \mathrm{rk}(\E) \deg(\E') - \mathrm{rk}(\E')\deg(\E)  \\
 & = &  \hspace{-6.7cm}\delta(\E', \E)  \\
 & \geq &  \hspace{-3.2cm} \mathrm{rk}(\E')(\mathrm{rk}(\E)-\mathrm{rk}(\E'))(2g-2)
\end{array}   $$
(with strict inequality in the stable case), we have $ H^1(X, (\E/\E')^{\vee} \otimes \E')=0$ by the previous lemma. 
Thus $\mathrm{Ext}^{1}(\E/\E', \E')=0$ and $\E \cong \E' \oplus \E/\E'.$ \end{proof}

\end{proposition}

This implies immediately the following. 

\begin{corollary} Let $\E$ be a $2$-bundle and $\Line \hookrightarrow \E$ a line subbundle with $\delta(\Line,\E) > 2g-2$, then $\E \cong \Line  \oplus \E/\Line. \hfill \square$

\end{corollary}

\begin{lemma} \label{lemma4.9} Let $\Line, \Line_1, \ldots, \Line_n$ be invertible sheaves. If 
$$\iota : \Line \longrightarrow \Line_1 \oplus \cdots \oplus \Line_n$$
is a non-zero morphism, then $\deg(\Line) \leq \deg(\Line_j)$ for some $j \in \{1, \ldots , n\}.$

\begin{proof} Write $\iota = \iota_1 \oplus \cdots \oplus \iota_n$, as $\iota$ is non-zero, there is some $j \in \{1, \ldots , n\}$ such that $\iota_j : \Line \rightarrow \Line_j$ is non-zero and thus injective. Therefore $\deg(\Line) \leq \deg(\Line_j).$ 
\end{proof}
\end{lemma}

\begin{proposition} Let $\Line_1, \ldots, \Line_n, \Line_{1}^{'}, \ldots, \Line_{m}^{'}$ be  invertible sheaves and
\[\iota : \Line_1 \oplus \cdots \oplus \Line_n  \rightarrow \Line_{1}^{'} \oplus \cdots \oplus \Line_{m}^{'}\]
be an injective morphism. Then there is an injection $ \{1, \ldots,n\} \rightarrow \{1. \ldots,m\}, i \mapsto j_i$, such that $\deg(\Line_i) \leq \deg(\Line_{j_i}^{'}).$

\begin{proof} Let  $d_i = \deg(\Line_i)$ and $d_{j}^{'} = \deg(\Line_{j}^{'})$ for  $i=1, \ldots,n$ and $j=1, \ldots, m.$ Suppose that $d_1 \leq d_2 \leq \cdots \leq d_n$ and $d_{1}^{'} \leq d_{2}'  \leq \cdots \leq d_{m}^{'}.$ 
We proceed by induction on $n$. The case $n=1$ follows from the previous lemma.  If $d_j' < d_1$ for all $j \in \{1, \ldots,m\}$, then $\mathrm{Hom}(\Line_i, \Line_j')= 0 $ for all $j=1,\ldots,m$, and $\iota$ is not injective. Let $j_0 = \min \{j| d_j' > d_1\}$, hence $\mathrm{Hom}(\Line_1 \oplus \cdots \oplus \Line_n, \Line_{j_0}') \neq 0$ and $m-j_0 \geq n$. Consider the commutative diagram with exact rows and columns 
$$\xymatrix@R5pt@C7pt{ && 0 \ar[dd] && 0 \ar[dd] && 0 \ar[dd] && \\ 
&& && && && \\
0 \ar[rr] && \iota^{*}\Line_{j_0}^{'} \ar[dd]^{\iota '}\ar[rr]^{\hspace{-0.8cm} f_1} && \Line_{1} \oplus \cdots \oplus \Line_{n} \ar[dd]^{\iota}\ar[rr]^{\hspace{0.4cm} f_2} && Q \ar@{-->}[dd]^{\varphi}\ar[rr] && 0 \   \ \textbf{(A)} \\
&& && && && \\
0 \ar[rr] && \Line_{j_0}^{'} \ar[rr]^{\hspace{-0.8cm} g_1} && \Line_{1}^{'} \oplus \cdots \oplus \Line_{m}^{'} \ar[rr]^{ \hspace{0.4cm} g_2} && \displaystyle\bigoplus_{\substack{ j=1 \\ j \neq j_0}}^{m}\Line_{j}' \ar[rr] && 0 \   \ \textbf{(B)}
}$$
where $Q:= (\Line_{1} \oplus \cdots \oplus \Line_{n}) / \iota^{*}\Line_{j_0}^{'}.$

Next we prove the existence and injectivity of $\varphi$. First observe that, since the diagram commutes, $g_2 \circ \iota \circ f_1 = 0$. The universal property of cokernels implies the existence of $\varphi$. It is enough to prove the injectivity on stalks (\cite{biblia} ex II.1.2). Let $x \in X$ and suppose $\varphi_x(q)=0$ for some  $q \in Q_x.$ Since $(f_2)_x$ is surjective, $q = (f_2)_x(a)$ for some  $a \in (\Line_{1} \oplus \cdots \oplus \Line_{n})_x,$ and we conclude successively 
\begin{itemize}
\item  $(g_2)_x((\iota)_x(a))= \varphi_x(q)= 0;$  
\item  $ (\iota)_x(a) \in \ker((g_2)_x) = \mathrm{im}((g_1)_x)  = (\Line_{j_0}^{'})_x;  $
\item  $ a \in( \Line_{1} \oplus \cdots \oplus \Line_{n})_x \text{ and } (\iota)_x(a) \in (\Line_{j_0}^{'})_x; $
\item  $ a \in (\iota^{*}\Line_{j_0}^{'})_x; $ 
\item  $q=(g_2)_x(a)=0.$
\end{itemize} 
This shows that $\varphi$ is injective. The \textit{snake lemma} implies that the following diagram is commutative with exact rows and columns

$$\xymatrix@R5pt@C7pt{ && 0 \ar[dd] && 0 \ar[dd] && 0 \ar[dd] && \\ 
&& && && && \\
0 \ar[rr] && \iota^{*}\Line_{j_0}^{'} \ar[dd]^{\iota '}\ar[rr]^{\hspace{-0.8cm} f_1} && \Line_{1} \oplus \cdots \oplus \Line_{n} \ar[dd]^{\iota}\ar[rr]^{\hspace{0.4cm} f_2} && Q \ar[dd]^{\varphi}\ar[rr] && 0 \   \ \textbf{(A)} \\
&& && && && \\
0 \ar[rr] && \Line_{j_0}^{'} \ar[dd] \ar[rr]^{\hspace{-0.8cm} g_1} && \Line_{1}^{'} \oplus \cdots \oplus \Line_{m}^{'} \ar[dd] \ar[rr]^{ \hspace{0.4cm} g_2} && \displaystyle\bigoplus_{\substack{ j=1 \\ j \neq j_0}}^{m}\Line_{j}' \ar[dd]\ar[rr] && 0 \   \ \textbf{(B)} \\
&& && && && \\
0 \ar[rr]&&\mathrm{coker}(\iota')\ar[dd] \ar[rr]&&\mathrm{coker}(\iota)\ar[dd] \ar[rr]&&\mathrm{coker}(\varphi) \ar[dd] \ar[rr]&& 0\\
&& && && && \\
&& 0 && 0 && 0 && \\
}$$

Since \textbf{(B)} splits, there is $\pi_1 : \Line_{1}^{'} \oplus \cdots \oplus \Line_{m}^{'} \rightarrow \Line_{j_0}^{'}$ such that $\pi_1 \circ g_1 = \mathrm{id}_{\Line_{j_0}^{'}}.$ 
The above diagram implies $\mathrm{im}(\pi_1 \circ \iota) \subset \mathrm{im}(\iota ').$ Indeed, (on stalks) if $y \in \mathrm{im}(\pi_1 \circ \iota) \setminus \mathrm{im}(\iota')$, then $y$ goes to non-zero element in $\mathrm{coker}(\iota).$ Hence $g_1(y) \not\in \mathrm{im}(\iota)$ which implies that $y \not\in \mathrm{im}(\pi_1 \circ  \iota)$ since $\pi_1 \circ g_1 = \mathrm{id}_{\Line_{j_0}^{'}}.$     

Define 
$$\pi_2 : \Line_{1} \oplus \cdots \oplus \Line_{n} \longrightarrow \iota^{*}\Line_{j_0}^{'}$$
by $\pi_2 = (\iota ')^{-1}  \circ \pi_1 \circ \iota$. Since $\iota '$ is injective,  $\pi_2$ is well defined. Moreover, $\pi_2 \circ f_1 = \mathrm{id}_{\iota^{*}\Line_{j_0}^{'}}$, thus \textbf{(A)} also is a split sequence and 
$$\Line_{1} \oplus \cdots \oplus \Line_{n} = \iota^{*}\Line_{j_0}^{'} \oplus Q.$$
Therefore $\iota^{*}\Line_{j_0}^{'} = \Line_{i}$ for some $i \in \{1, \ldots, n\}.$ By Lemma \ref{lemma4.9},  $d_i \leq d_{j_0}'.$ Thus the claim follows from the induction hypothesis.
\end{proof}

\end{proposition}

\begin{corollary} Let $\E = \Line_1 \oplus \cdots \oplus \Line_n$ be an $n$-bundle where $\Line_1, \ldots, \Line_n$ are line bundles and $\deg(\Line_1) \leq \cdots \leq \deg(\Line_n).$ Then 
$$\delta_k(\E) \geq n \sum_{i=0}^{k-1} \deg(\Line_{r-i}) - k \sum_{i=1}^{r} \deg(\Line_i)$$
for every $k=1, \ldots, n-1.$ In particular, $\delta_k(\E) \geq 0$ for $k=1, \ldots, n-1.$ \qed

\end{corollary}

\begin{corollary} If $r=2$ and $k=1$ in the last corollary, then equality holds, i.e.
\begin{center}
$\delta_1(\E) = \deg(\Line_2) - \deg(\Line_1).$ \qed  
\end{center} 
\end{corollary}

We use  the $\delta-$invariants to investigate vertices in the graph of an unramified Hecke operator which are connected by edge, i.e.\ to investigate sequences of the form
$$0 \longrightarrow \E' \longrightarrow \E \longrightarrow \mathcal{K}_{x}^{\oplus r} \longrightarrow 0$$ 
for $n-$bundles $\E$ and $\E'$,  $r=1, \ldots, n-1$, and $x \in X$ a closed point. 

If $\E'' \hookrightarrow \E$ is a subbundle, then we say that it \textit{lifts to} $\E'$ if there exists a morphism $\E'' \rightarrow \E',$ such that the diagram 
$$\xymatrix@R5pt@C7pt{ & & \E'' \ar[dd] \ar[ddll]\\  
& & \\
\E' \ar[rr] & & \E
}$$
commutes. In this case, $\E'' \rightarrow \E'$ is indeed a subbundle, since otherwise it would extend nontrivially to a subbundle $\overline{\E''} \rightarrow \E' \hookrightarrow \E$ and would contradict the hypothesis that $\E''$ is a subbundle of $\E.$ By exactness of the above sequence, a subbundle $\E'' \rightarrow \E$ lifts to $\E',$ if and only if the image of $\E''$ in $\mathcal{K}_{x}^{\oplus r} $ is $0$. 

\begin{theorem} \label{2maintheorem}Let $\E$ be an $n-$bundle. If $\E'$ is a neighbor of $\E$ in $\mathcal{G}_{x,r}$, then 
$$\delta_k(\E') \in \big\{ \delta_{k}(\E) -k |x|(n-r), \delta_{k}(\E) -k |x|(n-r) + n, \ldots,  \delta_{k}(\E)+  k |x|r\big\}$$
for $k=1, \ldots, n-1$. 

\begin{proof} Let $d=\deg(\E)$ and $d'=\deg(\E').$ By the exactness of the sequence
 $$0 \longrightarrow \E' \longrightarrow \E \longrightarrow \mathcal{K}_{x}^{\oplus r}  \longrightarrow 0,$$ 
we have $d' = d-r |x|$. Let $\F \hookrightarrow \E$  and $\F' \hookrightarrow \E'$ be  maximal $k$-subbundles. Note that every subbundle of $\E'$ is a locally free subsheaf of $\E$ and thus extends to a subbundle of $\E$. Let $\overline{\F'}$ be the subbundle of $\E$ that extends $\F'.$ We know  that $\deg(\overline{\F'}) \geq \deg(\F').$ Thus
\begin{align*}
\delta_{k}(\E') & =   n \deg(\F') - k d'  \\ 
  &\leq  n \deg(\overline{\F'}) -k(d-r |x|) \\ 
  &=  \delta_{k}(\overline{\F'}, \E) +  k |x| r \\ 
  &\leq  \delta_{k}(\E) +  k |x|r.
\end{align*}

If $\F \hookrightarrow \E$ lifts to $\E',$ then $\F \hookrightarrow \E'$ is a maximal $k$-subbundle. 
Thus we have equalities in the above estimation, i.e.\ $\delta_{k}(\E') = \delta_{k}(\E) +  k |x|r.$

Let $\Line(-x)$ be the ideal sheaf of $\{x\}.$  For every $k$-subbundle $\F \hookrightarrow \E$, we may think of  $\Line(-x) \otimes \F$ as a subsheaf of $\F.$ 

If $\F \hookrightarrow \E$ does not lift to $\E'$, then $\Line(-x) \otimes \F \subseteq \F \hookrightarrow \E$ lifts to a subsheaf of $\E'$ since $\Line(-x)= \ker(\mathcal{O}_X \rightarrow \mathcal{K}_x).$ Hence
$$\xymatrix@R5pt@C7pt{ \Line(-x) \otimes \F \ar[dd] & \subseteq & \F \ar[dd] \\
& & \\
\E' \ar[rr] & & \E  
}$$     
Therefore, 
\begin{align*}
\delta_k(\E') & =  \delta_{k}(\F',\E') \\ 
 & \geq  \delta_{k}(\overline{\Line(-x) \otimes \F}, \E') \\ 
 & =  n \deg(\overline{\Line(-x) \otimes \F}) - k d' \\ 
 & \geq  n(\deg(\F) - k |x|) - k(d - r |x|) \\
 & =  \delta_{k}(\F,\E) - n k |x| +  k |x| r\\
 & =  \delta_{k}(\E) - k |x|(n-r).
\end{align*}

Furthermore, 
$$\delta_{k}(\E') \equiv - k d' \equiv (d - r |x|)k \equiv \delta_{k}(\E) +  k |x| r \equiv \delta_{k}(\E) - (n-r)k |x| \   \  (\mathrm{mod}\  \ n).$$
This concludes the proof of the theorem. \end{proof}
\end{theorem}

\begin{corollary} Let $\E$ be an $n-$bundle and $\F \hookrightarrow \E$ a maximal $k$-subbundle. If $\E'$ is a neighbor of $\E$ in $\mathcal{G}_{x,r}$ such that $\F \hookrightarrow \E$ lifts to $\E'$, then
$$\delta_k(\E') = \delta_{k}(\E) +  k |x|r.$$

\begin{proof} This follows from last theorem's proof.
\end{proof} 
\end{corollary}

\section{Graphs for the projective line}\label{section5}

With the theory from the previous sections, we are able to calculate all graphs for unramified Hecke operators over a rational function field. For the rest of the section, we fix $F = \mathbb{F}_q(T),$ i.e.\ $X$ is the projective line over $\mathbb{F}_q.$ We intend to determine the graphs of $\Phi_{x,r}$ for every $n \in \Z_{> 0}$, $r=1, \ldots, n$ and $x \in \mathbb{P}^1$ a closed point.  

First of all, we know by a theorem due to Grothendieck\footnote{actually it was already known by Dedekind and Weber, see \cite{dedekind-weber}.} (see \cite{Gortz}, Theorem 11.51) that every rank-$n$ vector bundle over $\mathbb{P}^1$ is isomorphic to
$$\mathcal{O}_{\mathbb{P}^1}(d_1) \oplus \cdots \oplus \mathcal{O}_{\mathbb{P}^1}(d_n)$$
for some $d_1 \geq \cdots \geq d_n.$ By Proposition \ref{propedge} and Grothendieck's theorem, the graphs for the zero element are given by a subset of the  standard $n-$lattice $\Z^n$ consisting of vertices without edges and the graphs for identity $1$ in $\mathcal{H}_K$ is a subset of the standard $n-$lattice $\Z^n$ consisting of vertices with loops of multiplicity one.

In certain cases, for a better reading we will denote by $(d_1, \ldots, d_n)$ the rank-$n$ vector bundle $\mathcal{O}_{\mathbb{P}^1}(d_1) \oplus \cdots \oplus \mathcal{O}_{\mathbb{P}^1}(d_n).$

\begin{figure}[h]
  \beginpgfgraphicnamed{tikz/graph001}
  \begin{tikzpicture}[>=latex, scale=2]
        \vertex[circle,fill,label={below:$\mathcal{O}_{\mathbb{P}^1}(-2)$}](00) at (-2,0) {};
        \vertex[circle,fill,label={below:$\mathcal{O}_{\mathbb{P}^1}(-1)$}](11) at (-1,0) {};
        \vertex[circle,fill,label={below:$\mathcal{O}_{\mathbb{P}^1}$}](10) at (0,0) {};
        \vertex[circle,fill,label={below:$\mathcal{O}_{\mathbb{P}^1}(1)$}](22) at (1,0) {};
        \vertex[circle,fill,label={below:$\mathcal{O}_{\mathbb{P}^1}(2)$}](21) at (2,0) {};
        \draw (-2.5,0) circle (0.015cm);
\fill  (-2.5,0) circle (0.015cm);
  \draw (-2.7,0) circle (0.015cm);
\fill  (-2.7,0) circle (0.015cm);  
\draw (-2.9,0) circle (0.015cm);
\fill  (-2.9,0) circle (0.015cm);    
       \draw (2.5,0) circle (0.015cm);
\fill  (2.5,0) circle (0.015cm);
  \draw (2.7,0) circle (0.015cm);
\fill  (2.7,0) circle (0.015cm);  
\draw (2.9,0) circle (0.015cm);
\fill  (2.9,0) circle (0.015cm);  
        \end{tikzpicture}
 \endpgfgraphicnamed
 \caption{The graph of zero element in $\mathcal{H}_K$ for $n=1$. }
  \label{graph001}
\end{figure}
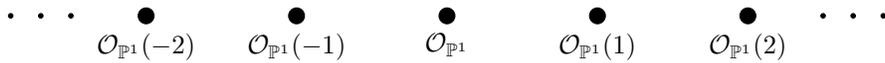
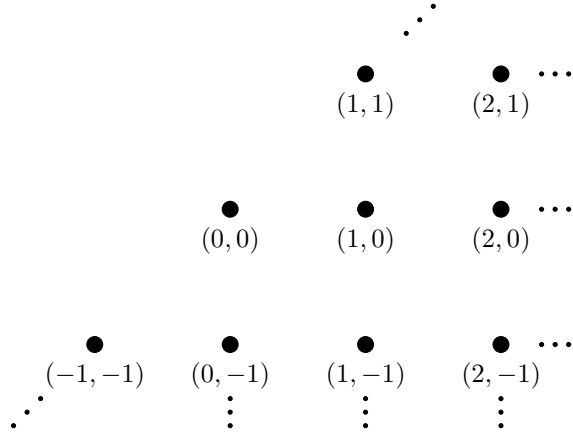
\begin{figure}[h]
  \beginpgfgraphicnamed{tikz/graph02}
  \begin{tikzpicture}[>=latex, scale=1.8]
        \vertex[circle,fill,label={below:$(0,0)$}](00) at (0,0) {};
        \vertex[circle,fill,label={below:$(1,1)$}](11) at (1,1) {};
        \vertex[circle,fill,label={below:$(-1,-1)$}](10) at (-1,-1) {};
        \vertex[circle,fill,label={below:$(1,0) $}](22) at (1,0) {};
        \vertex[circle,fill,label={below:$(2,0)$}](21) at (2,0) {};
         \vertex[circle,fill,label={below:$(0,-1)$}](21) at (0,-1) {};
          \vertex[circle,fill,label={below:$(1,-1)$}](21) at (1,-1) {};
           \vertex[circle,fill,label={below:$(2,-1)$}](21) at (2,-1) {};
            \vertex[circle,fill,label={below:$(2,1)$}](21) at (2,1) {}; 
 \draw (-1.4,-1.4) circle (0.015cm);
\fill  (-1.4,-1.4) circle (0.015cm);
  \draw (-1.5,-1.5) circle (0.015cm);
\fill  (-1.5,-1.5) circle (0.015cm);  
\draw (-1.6,-1.6) circle (0.015cm);
\fill  (-1.6,-1.6) circle (0.015cm);  

 \draw (1.3,1.3) circle (0.015cm);
\fill  (1.3,1.3) circle (0.015cm);
  \draw (1.4,1.4) circle (0.015cm);
\fill  (1.4,1.4) circle (0.015cm);  
\draw (1.5,1.5) circle (0.015cm);
\fill  (1.5,1.5) circle (0.015cm);  
 
 \draw (0,-1.4) circle (0.015cm);
\fill  (0,-1.4) circle (0.015cm);
  \draw (0,-1.5) circle (0.015cm);
\fill  (0,-1.5) circle (0.015cm);  
\draw (0,-1.6) circle (0.015cm);
\fill  (0,-1.6) circle (0.015cm); 

 \draw (1,-1.4) circle (0.015cm);
\fill  (1,-1.4) circle (0.015cm);
  \draw (1,-1.5) circle (0.015cm);
\fill  (1,-1.5) circle (0.015cm);  
\draw (1,-1.6) circle (0.015cm);
\fill  (1,-1.6) circle (0.015cm); 

 \draw (2,-1.4) circle (0.015cm);
\fill  (2,-1.4) circle (0.015cm);
  \draw (2,-1.5) circle (0.015cm);
\fill  (2,-1.5) circle (0.015cm);  
\draw (2,-1.6) circle (0.015cm);
\fill  (2,-1.6) circle (0.015cm); 

 \draw (2.3,-1) circle (0.015cm);
\fill  (2.3,-1) circle (0.015cm);
  \draw (2.4,-1) circle (0.015cm);
\fill  (2.4,-1) circle (0.015cm);  
\draw (2.5,-1) circle (0.015cm);
\fill  (2.5,-1) circle (0.015cm); 

 \draw (2.3,0) circle (0.015cm);
\fill  (2.3,0) circle (0.015cm);
  \draw (2.4,0) circle (0.015cm);
\fill  (2.4,0) circle (0.015cm);  
\draw (2.5,0) circle (0.015cm);
\fill  (2.5,0) circle (0.015cm); 

 \draw (2.3,1) circle (0.015cm);
\fill  (2.3,1) circle (0.015cm);
  \draw (2.4,1) circle (0.015cm);
\fill  (2.4,1) circle (0.015cm);  
\draw (2.5,1) circle (0.015cm);
\fill  (2.5,1) circle (0.015cm); 
        \end{tikzpicture}
 \endpgfgraphicnamed
\caption{The graph of zero element in $\mathcal{H}_K$ for $n=2$.} 
 \label{graph02}
\end{figure}
\begin{figure}[h]
  \beginpgfgraphicnamed{tikz/graph11}
  \begin{tikzpicture}[>=latex, scale=2]
        \vertex[circle,fill,label={below:${\tiny\mathcal{O}_{\mathbb{P}^1}}$}](00) at (0,0) {};
  \draw[->-=0.5] (0,0.25) circle  (0.25cm) node at (0,0.6) {\normalsize $1$} ;
  
     \vertex[circle,fill,label={below:${\tiny \mathcal{O}_{\mathbb{P}^1}(-2)}$}](00) at (-2,0) {};
  \draw[->-=0.5] (-2,0.25) circle  (0.25cm) node at (-2,0.6) {\normalsize $1$}; 
  
     \vertex[circle,fill,label={below:${\tiny\mathcal{O}_{\mathbb{P}^1}(-1)}$}](00) at (-1,0) {};
  \draw[->-=0.5] (-1,0.25) circle  (0.25cm) node at (-1,0.6) {\normalsize $1$} ;
  
     \vertex[circle,fill,label={below:${\tiny\mathcal{O}_{ \mathbb{P}^1}(1)}$}](00) at (1,0) {};
  \draw[->-=0.5] (1,0.25) circle  (0.25cm) node at (1,0.6) {\normalsize $1$} ; 
  
     \vertex[circle,fill,label={below:${\tiny\mathcal{O}_{ \mathbb{P}^1}(2)}$}](00) at (2,0) {};
  \draw[->-=0.5] (2,0.25) circle  (0.25cm) node at (2,0.6) {\normalsize $1$} ;

 \draw (2.4,0) circle (0.015cm);
\fill  (2.4,0) circle (0.015cm);
  \draw (2.5,0) circle (0.015cm);
\fill  (2.5,0) circle (0.015cm);  
\draw (2.6,0) circle (0.015cm);
\fill  (2.6,0) circle (0.015cm); 

 \draw (-2.4,0) circle (0.015cm);
\fill  (-2.4,0) circle (0.015cm);
  \draw (-2.5,0) circle (0.015cm);
\fill  (-2.5,0) circle (0.015cm);  
\draw (-2.6,0) circle (0.015cm);
\fill  (-2.6,0) circle (0.015cm); 
   \end{tikzpicture}
 \endpgfgraphicnamed
 \caption{The graph of the identity in $\mathcal{H}_K$ for $n=1$.}
\label{graph11}
\end{figure}
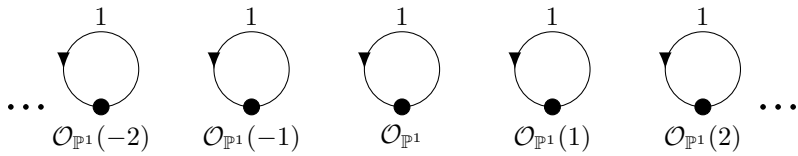
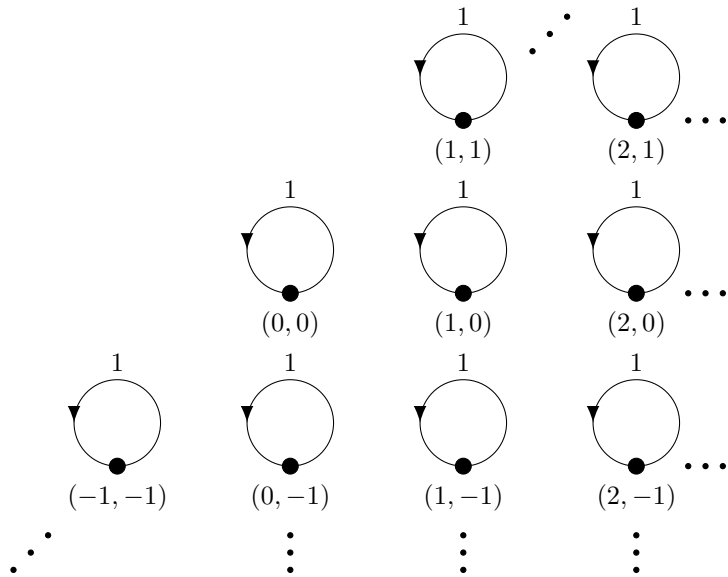
\begin{figure}[h]
  \beginpgfgraphicnamed{tikz/graph12}
  \begin{tikzpicture}[>=latex, scale=2.3]
        \vertex[circle,fill,label={below:$(0,0)$}](00) at (0,0) {};
         \draw[->-=0.5] (0,0.25) circle  (0.25cm) node at (0,0.6) {\normalsize $1$} ;
   
        \vertex[circle,fill,label={below:$(1,1)$}](11) at (1,1) {};
         \draw[->-=0.5] (1,1.25) circle  (0.25cm) node at (1,1.6) {\normalsize $1$} ;
         
        \vertex[circle,fill,label={below:$(-1,-1)$}](10) at (-1,-1) {};
         \draw[->-=0.5] (-1,-0.75) circle  (0.25cm) node at (-1,-0.4) {\normalsize $1$} ;
         
        \vertex[circle,fill,label={below:$(1,0) $}](22) at (1,0) {};
         \draw[->-=0.5] (1,0.25) circle  (0.25cm) node at (1,0.6) {\normalsize $1$} ;

        \vertex[circle,fill,label={below:$(2,0)$}](21) at (2,0) {};
        \draw[->-=0.5] (2,0.25) circle  (0.25cm) node at (2,0.6) {\normalsize $1$} ; 
         
         \vertex[circle,fill,label={below:$(0,-1)$}](21) at (0,-1) {};
           \draw[->-=0.5] (0,-0.75) circle  (0.25cm) node at (0,-0.4) {\normalsize $1$} ;
           
          \vertex[circle,fill,label={below:$(1,-1)$}](21) at (1,-1) {};
 \draw[->-=0.5] (1,-0.75) circle  (0.25cm) node at (1,-0.4) {\normalsize $1$} ;
            
           \vertex[circle,fill,label={below:$(2,-1)$}](21) at (2,-1) {};
 \draw[->-=0.5] (2,-0.75) circle  (0.25cm) node at (2,-0.4) {\normalsize $1$} ;
 
            \vertex[circle,fill,label={below:$(2,1)$}](21) at (2,1) {}; 
 \draw[->-=0.5] (2,1.25) circle  (0.25cm) node at (2,1.6) {\normalsize $1$} ; 
 
 \draw (-1.4,-1.4) circle (0.015cm);
\fill  (-1.4,-1.4) circle (0.015cm);
  \draw (-1.5,-1.5) circle (0.015cm);
\fill  (-1.5,-1.5) circle (0.015cm);  
\draw (-1.6,-1.6) circle (0.015cm);
\fill  (-1.6,-1.6) circle (0.015cm);  

 \draw (1.4,1.4) circle (0.015cm);
\fill  (1.4,1.4) circle (0.015cm);
  \draw (1.5,1.5) circle (0.015cm);
\fill  (1.5,1.5) circle (0.015cm);  
\draw (1.6,1.6) circle (0.015cm);
\fill  (1.6,1.6) circle (0.015cm);  
 
 \draw (0,-1.4) circle (0.015cm);
\fill  (0,-1.4) circle (0.015cm);
  \draw (0,-1.5) circle (0.015cm);
\fill  (0,-1.5) circle (0.015cm);  
\draw (0,-1.6) circle (0.015cm);
\fill  (0,-1.6) circle (0.015cm); 

 \draw (1,-1.4) circle (0.015cm);
\fill  (1,-1.4) circle (0.015cm);
  \draw (1,-1.5) circle (0.015cm);
\fill  (1,-1.5) circle (0.015cm);  
\draw (1,-1.6) circle (0.015cm);
\fill  (1,-1.6) circle (0.015cm); 

 \draw (2,-1.4) circle (0.015cm);
\fill  (2,-1.4) circle (0.015cm);
  \draw (2,-1.5) circle (0.015cm);
\fill  (2,-1.5) circle (0.015cm);  
\draw (2,-1.6) circle (0.015cm);
\fill  (2,-1.6) circle (0.015cm); 

 \draw (2.3,-1) circle (0.015cm);
\fill  (2.3,-1) circle (0.015cm);
  \draw (2.4,-1) circle (0.015cm);
\fill  (2.4,-1) circle (0.015cm);  
\draw (2.5,-1) circle (0.015cm);
\fill  (2.5,-1) circle (0.015cm); 

 \draw (2.3,0) circle (0.015cm);
\fill  (2.3,0) circle (0.015cm);
  \draw (2.4,0) circle (0.015cm);
\fill  (2.4,0) circle (0.015cm);  
\draw (2.5,0) circle (0.015cm);
\fill  (2.5,0) circle (0.015cm); 

 \draw (2.3,1) circle (0.015cm);
\fill  (2.3,1) circle (0.015cm);
  \draw (2.4,1) circle (0.015cm);
\fill  (2.4,1) circle (0.015cm);  
\draw (2.5,1) circle (0.015cm);
\fill  (2.5,1) circle (0.015cm); 
        \end{tikzpicture}
 \endpgfgraphicnamed
\caption{The graph of the identity in $\mathcal{H}_K$ for $n=2$.}
 \label{graph12}
\end{figure}

We draw these two cases below for $n=1,2$, see figures $\ref{graph001}$, $\ref{graph02}$, $\ref{graph11}$ and $\ref{graph12}$ respectively.  In certain cases, for a better reading we will denote by $(d_1, \ldots, d_n)$ the rank $n$-vector bundle $\mathcal{O}_{\mathbb{P}^1}(d_1) \oplus \cdots \oplus \mathcal{O}_{\mathbb{P}^1}(d_n).$\\

\newpage


We continue with the calculation of the graph of $\Phi_{x,1}$, for every closed point $x$ in $\mathbb{P}^1$ of degree $|x|$ and $n=1.$ By the correspondence of Theorem \ref{theoremcorrespondence}, a neighborhood of $\mathcal{O}_{\mathbb{P}^1}(d)$ is given in terms of exact sequences
$$0 \longrightarrow \E' \longrightarrow \mathcal{O}_{\mathbb{P}^1}(d) \longrightarrow \mathcal{K}_x \longrightarrow 0.$$ 
By Grothendieck's classification of vector bundles on $\mathbb{P}^{1}$, $\E'$ is completely determined by its degree. By additivity of the degree, $\deg \E' = d - |x|,$
thus $\E' = \mathcal{O}_{\mathbb{P}^1}(d-|x|)$ and 
$$m_{x,1}(\mathcal{O}_{\mathbb{P}^1}(d),\mathcal{O}_{\mathbb{P}^1}(d-|x|)) =1$$ 
(Corollary  \ref{corollaryofmaintheorem}).  The graph is illustrated as in the figure \ref{graph,n=1}.

\begin{figure}[h]
  \beginpgfgraphicnamed{tikz/graph,n=1}
  \begin{tikzpicture}[>=latex, scale=2]
        \vertex[circle,fill,label={below:$\mathcal{O}_{\mathbb{P}^1}(d)$}](00) at (-2,0) {};
        \vertex[circle,fill,label={below:$\mathcal{O}_{\mathbb{P}^1}(d-|x|)$}](11) at (0,0) {};
        \vertex[circle,fill,label={below:$\mathcal{O}_{\mathbb{P}^1}(d-2|x|)$}](10) at (2,0) {};
        \vertex[circle,fill,label={below:$\mathcal{O}_{\mathbb{P}^1}(d-1)$}](22) at (-2,-0.6) {};
        \vertex[circle,fill,label={below:$\mathcal{O}_{\mathbb{P}^1}(d-1-|x|)$}](21) at (0,-0.6) {};
        \vertex[circle,fill,label={below:$\mathcal{O}_{\mathbb{P}^1}(d-1-2|x|)$}](20) at (2,-0.6) {};
        
                \vertex[circle,fill,label={below:$\mathcal{O}_{\mathbb{P}^1}(d-|x|+1)$}](3) at (-2,-1.7) {};
        \vertex[circle,fill,label={below:$\mathcal{O}_{\mathbb{P}^1}(d-2|x|+1)$}](4) at (0,-1.7) {};
        \vertex[circle,fill,label={below:{$\tiny \mathcal{O}_{\mathbb{P}^1}(d-3|x|+1)$}}](5) at (2,-1.7) {};
    \path[-,font=\scriptsize]
    (00) edge[->-=0.8] node[pos=0.2,auto,black] {\normalsize $1$} (11)
    (11) edge[->-=0.8] node[pos=0.2,auto,black] {\normalsize $1$} (10)
    (22) edge[->-=0.8] node[pos=0.2,auto,black] {\normalsize $1$} (21)
    (21) edge[->-=0.8] node[pos=0.2,auto,black] {\normalsize $1$} (20)
    
    (3) edge[->-=0.8] node[pos=0.2,auto,black] {\normalsize $1$} (4)
    (4) edge[->-=0.8] node[pos=0.2,auto,black] {\normalsize $1$} (5)
    
    (10) edge[-] node[pos=0.2,auto,black] {} (2.3,0)
     (20) edge[-] node[pos=0.2,auto,black] {} (2.3,-0.6)
      (-2.3,0) edge[-] node[pos=0.2,auto,black] {} (00)
       (-2.3,-0.6) edge[-] node[pos=0.2,auto,black] {} (22)
       
(5) edge[-] node[pos=0.2,auto,black] {} (2.3,-1.7)       
(-2.3,-1.7) edge[-] node[pos=0.2,auto,black] {} (3)       
       
       ;
        \draw (2.4,0) circle (0.015cm);
\fill  (2.4,0) circle (0.015cm);
  \draw (2.5,0) circle (0.015cm);
\fill  (2.5,0) circle (0.015cm);  
\draw (2.6,0) circle (0.015cm);
\fill  (2.6,0) circle (0.015cm);

        \draw (2.4,-1.7) circle (0.015cm);
\fill  (2.4,-1.7) circle (0.015cm);
  \draw (2.5,-1.7) circle (0.015cm);
\fill  (2.5,-1.7) circle (0.015cm);  
\draw (2.6,-1.7) circle (0.015cm);
\fill  (2.6,-1.7) circle (0.015cm);

        \draw (-2.4,-1.7) circle (0.015cm);
\fill  (-2.4,-1.7) circle (0.015cm);
  \draw (-2.5,-1.7) circle (0.015cm);
\fill  (-2.5,-1.7) circle (0.015cm);  
\draw (-2.6,-1.7) circle (0.015cm);
\fill  (-2.6,-1.7) circle (0.015cm);

        \draw (2.4,-0.6) circle (0.015cm);
\fill  (2.4,-0.6) circle (0.015cm);
  \draw (2.5,-0.6) circle (0.015cm);
\fill  (2.5,-0.6) circle (0.015cm);  
\draw (2.6,-0.6) circle (0.015cm);
\fill  (2.6,-0.6) circle (0.015cm);

        \draw (-2.4,0) circle (0.015cm);
\fill  (-2.4,0) circle (0.015cm);
  \draw (-2.5,0) circle (0.015cm);
\fill  (-2.5,0) circle (0.015cm);  
\draw (-2.6,0) circle (0.015cm);
\fill  (-2.6,0) circle (0.015cm);

        \draw (-2.4,-0.6) circle (0.015cm);
\fill  (-2.4,-0.6) circle (0.015cm);
  \draw (-2.5,-0.6) circle (0.015cm);
\fill  (-2.5,-0.6) circle (0.015cm);  
\draw (-2.6,-0.6) circle (0.015cm);
\fill  (-2.6,-0.6) circle (0.015cm);

        \draw (-1,-1) circle (0.015cm);
\fill  (-1,-1) circle (0.015cm);
  \draw (-1,-1.1) circle (0.015cm);
\fill  (-1,-1.1) circle (0.015cm);  
\draw (-1,-1.2) circle (0.015cm);
\fill  (-1,-1.2) circle (0.015cm);

        \draw (1,-1) circle (0.015cm);
\fill  (1,-1) circle (0.015cm);
  \draw (1,-1.1) circle (0.015cm);
\fill  (1,-1.1) circle (0.015cm);  
\draw (1,-1.2) circle (0.015cm);
\fill  (1,-1.2) circle (0.015cm); 

   \end{tikzpicture}
 \endpgfgraphicnamed
 \caption{The graph of $\Phi_{x,1}$ in $\mathcal{H}_K$ for $n=1$.}
 \label{graph,n=1}
\end{figure}
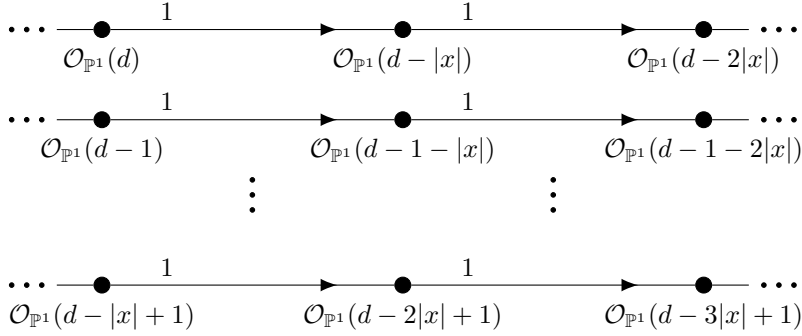


Our next aim is to explain an algorithm to calculate the graphs for $\Phi_{x,r}$ and $n>1.$ First we treat the case $|x|=1$ and afterwards we treat arbitrary degrees.

The closed points $y \in \mathbb{P}^1$ are in one-yo-one correspondence with irreducible and monic polynomials $\pi_y$ in $\mathbb{F}_q[T]$ and the infinite valuation. Moreover $ |y| = \deg \pi_y$ and $\pi_y$ is the uniformizer of $y.$ For $z=f(T)/g(T) \in  F \setminus \{0\} $ with $(f,g)=1,$
$$[F:\mathbb{F}_q(z)] = \mathrm{max}\{\deg(f),\deg(g)\}. $$
If $x$ is a degree one place, then $[F:\mathbb{F}_q(\pi_x)]=\deg(\pi_x)=1.$
Therefore we can suppose $T=\pi_x.$ Hence we can cover $\mathbb{P}^1$ by two open affine subschemes $U_0 := \mathrm{Spec} R^{+}$ and $U_1 := \mathrm{Spec} R^{-}$, where 
$$R^{+} := \mathbb{F}_q[\pi_x], \   \ R^{-} := \mathbb{F}_q[\pi_{x}^{-1}] \text{  and  } R^{\pm} := \mathbb{F}_q[\pi_x, \pi_{x}^{-1}].$$
We have  bijections 
$$
G(R^{+}) \backslash G(R^{\pm}) / G(R^{-})     \cong   \mathrm{Bun}_n \mathbb{P}^1 \cong G(F)\backslash G(\mathbb{A}) / K  $$ \vspace{-0.7cm}
$$g_x  \hspace{2cm}\longrightarrow  \hspace{1cm} (g_x, \mathrm{id})$$
where $(g_x, \mathrm{id})$ means the class of an adelic matrix which equals $g_x$ at place the $x$ and the identity matrix at all places $y \neq x$ (see \cite{Gortz} Chapter 11 for more details).

\begin{lemma}[\cite{Gortz} Lemma 11.50] \label{lemmaG}  Let $(\Z^n)_{+}$ be the  set of $\mathrm{\textbf{\emph{d}}} = (d_i)_i \in \Z^n$ with $d_1 \geq d_2 \geq \cdots \geq d_n.$ For each $\mathrm{\textbf{\emph{d}}} \in (\Z^n)_{+}$, let $\pi_{x}^{\mathrm{\textbf{\emph{d}}}} \in G(R^{\pm})$ be the diagonal matrix with entries $\pi_{x}^{d_1}, \ldots, \pi_{x}^{d_n}.$ Then the following map is bijective
$$\begin{array}{ccc}
(\Z^n)_{+} & \longrightarrow & G(R^{-}) \backslash G(R^{\pm}) / G(R^{+}). \\ \vspace{0.3cm}
\mathrm{\textbf{\emph{d}}}  & \longmapsto &  G(R^{-})  \pi_{x}^{\mathrm{\textbf{\emph{d}}}}  G(R^{+})
\end{array} $$

\end{lemma}

In the case of a degree one place $x$, we can determine the graph of  $\Phi_{x,r}$, for  $n > 1$ and $r \in \{1, \ldots, n\}$ by finding the representative $G(R^{-}) \pi_{x}^{\textbf{\emph{d}}'} G(R^{+})$ for the double cosets of $ \pi_{x}^{\textbf{\emph{d}}} \xi_w, $ cf. Lemma \ref{lemma1} and Theorem \ref{maintheorem}.


We shall use the symbol "$\sim$"  when two matrices represent the same class  in  $ G(R^{-}) \backslash G(R^{\pm}) / G(R^{+}).$

\begin{example} In the following, for $n=2$ i.e.\ rank two bundles, we describe the graph of $\Phi_{x,1}$ where $x$ a degree one place.  By Theorem \ref{maintheorem} and Lemma \ref{lemmaG}, we have to find the unique diagonal matrix $\pi_{x}^{\mathrm{\textbf{d}}}$ with $ \mathrm{\textbf{d}} \in  (\Z^n)_{+}$ that represents $g\xi_w$ in   $ G(R^{-}) \backslash G(R^{\pm}) / G(R^{+})$ where 
$$[g] = \left(\begin{small}
 \begin{array}{cc}
\pi_{x}^{d_1} &  \\ 
 & \pi_{x}^{d_2} 
\end{array} \end{small}\right)  $$
with $d_1 \geq d_2$ which corresponds to the vector bundle $\mathcal{O}_{\mathbb{P}^1}(d_1) \oplus \mathcal{O}_{\mathbb{P}^1}(d_2),$  and where 
$$\xi_w = \left( \begin{small} \begin{array}{cc}
\pi_{x} & b \\ 
 & 1
\end{array} \end{small}\right) \text{ with } w=[1:b] \in \mathrm{Gr}(1,2)(\kappa(x)) = \mathbb{P}^1(\mathbb{F}_q)$$
or
$$\xi_w = \left(\begin{small} \begin{array}{cc}
1 &  \\ 
 & \pi_{x}^{}
\end{array} \end{small}\right) \text{ with } w=[0:1]  \in \mathrm{Gr}(1,2)(\kappa(x)) = \mathbb{P}^1(\mathbb{F}_q). $$
For $d_1> d_2, w=[1:b] \text{ and }  b \neq 0$, we have
$$g \xi_w = \left(\begin{small}  \begin{array}{cc}
\pi_{x}^{d_1} &  \\ 
 & \pi_{x}^{d_2} 
\end{array}\end{small} \right) \left( \begin{array}{cc}
\pi_{x} & b \\ 
 & 1 
\end{array} \right) = \left( \begin{array}{cc}
\pi_{x}^{d_1+1} & b \pi_{x}^{d_1} \\ 
 & \pi_{x}^{d_2} 
\end{array} \right)  $$  
$$\sim \left( \begin{array}{cc}
 & 1 \\ 
1 & 
\end{array} \right)  \left( \begin{array}{cc}
\pi_{x}^{d_1+1} & b \pi_{x}^{d_1} \\ 
 & \pi_{x}^{d_2} 
\end{array} \right) \left( \begin{array}{cc}
 & 1 \\ 
1 & 
\end{array} \right) = 
\left( \begin{array}{cc}
\pi_{x}^{d_2} &  \\ 
b \pi_{x}^{d_1} & \pi_{x}^{d_1+1}
\end{array} \right) $$
$$\sim \left( \begin{array}{cc}
1 & -b^{-1} \pi_{x}^{d_2-d_1} \\ 
 & 1\
\end{array} \right) 
\left( \begin{array}{cc}
\pi_{x}^{d_2}&  \\ 
b \pi_{x}^{d_1} & \pi_{x}^{d_1+1}\
\end{array} \right) =\left( \begin{array}{cc}
 &  -b^{-1}\pi_{x}^{d_2+1}\\ 
b \pi_{x}^{d_1} & \pi_{x}^{d_1+1} \
\end{array} \right)  $$
$$\sim \left( \begin{array}{cc}
 & 1 \\ 
1 & 
\end{array} \right) \left( \begin{array}{cc}
 &  -b^{-1}\pi_{x}^{d_2+1}\\ 
b \pi_{x}^{d_1} & \pi_{x}^{d_1+1} \
\end{array} \right) = \left( \begin{array}{cc}
b \pi_{x}^{d_1} & \pi_{x}^{d_1+1} \\ 
 & -b^{-1}\pi_{x}^{d_2+1}\
\end{array} \right)  $$
$$ \sim \left( \begin{array}{cc}
b \pi_{x}^{d_1} & \pi_{x}^{d_1+1} \\ 
 & -b^{-1} \pi_{x}^{d_2+1}
\end{array} \right)
\left( \begin{array}{cc}
1 & -b^{-1} \pi_{x}^{} \\ 
 & 1
\end{array} \right) =  \left( \begin{array}{cc}
b \pi_{x}^{d_1} &  \\ 
 & -b^{-1}\pi_{x}^{d_2+1} 
\end{array} \right)  $$
$$\sim \left( \begin{array}{cc}
b\pi_{x}^{d_1} &  \\ 
 & -b^{-1} \pi_{x}^{d_2+1}
\end{array} \right)
\left( \begin{array}{cc}
b^{-1} &  \\ 
 & -b
\end{array} \right) = 
\left( \begin{array}{cc}
\pi_{x}^{d_1} &  \\ 
 & \pi_{x}^{d_2+1}
\end{array} \right) =: g_1.$$
Note that the number of classes $[g \xi_w]$ with $w=[1:b$] and $b \neq 0$ is $q-1.$

For $b=0,$ we have
$$g \xi_w = \left( \begin{array}{cc}
\pi_{x}^{d_1+1} &  \\ 
 & \pi_{x}^{d_2} 
\end{array} \right) =: g_2.$$
Thus $[g_2]$ is a neighbor of $[g]$ of at least multiplicity $1.$

If $d_1> d_2$ and $w =[0:1],$ then
$$g \xi_w = \left( \begin{array}{cc}
\pi_{x}^{d_1} &  \\ 
 & \pi_{x}^{d_2} 
\end{array} \right) 
\left( \begin{array}{cc}
1 &  \\ 
 & \pi_{x} 
\end{array} \right) = 
\left( \begin{array}{cc}
\pi_{x}^{d_1} &  \\ 
 & \pi_{x}^{d_2+1} 
\end{array} \right) = g_1$$
with multiplicity $1.$ We conclude that for $d_1 > d_2$, the neighborhood of $[g]$ is 
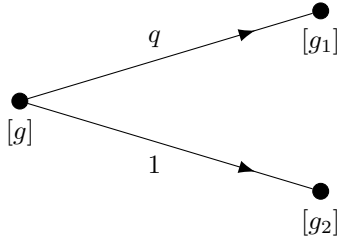
\[
  \beginpgfgraphicnamed{tikz/graphex1}
  \begin{tikzpicture}[>=latex, scale=2]
        \vertex[circle,fill,label={below:$[g]$}](00) at (0,0) {};
        \vertex[circle,fill,label={below:$[g_1]$}](11) at (2,0.6) {};
        \vertex[circle,fill,label={below:$[g_2]$}](10) at (2,-0.6) {};
  
            \path[-,font=\scriptsize]
    (00) edge[->-=0.8] node[pos=0.5,auto,black,swap] {\normalsize $1$} (10)
    (00) edge[->-=0.8] node[pos=0.5,auto,black] {\normalsize  $q$} (11)
    ;
   \end{tikzpicture}
 \endpgfgraphicnamed
\]

Observe that the sum of multiplicities in the vertex $[g]$ is $q+1 = \#\mathrm{Gr}(1,2)(\kappa(x))$, in accordance with Theorem \ref{maintheorem}.
If $d_1=d_2$, then $[g_1]=[g_2]$, and we have 

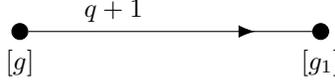
\[
  \beginpgfgraphicnamed{tikz/graphex1}
  \begin{tikzpicture}[>=latex, scale=2]
        \vertex[circle,fill,label={below:$[g]$}](00) at (0,0) {};
        \vertex[circle,fill,label={below:$[g_1]$}](11) at (2,0) {};
           \path[-,font=\scriptsize]
   (00) edge[->-=0.8] node[pos=0.3,auto,black] {\normalsize  $q+1$} (11)
    ;
   \end{tikzpicture}
 \endpgfgraphicnamed
\]

\end{example}

\begin{example} Next we illustrate the graph of $\Phi_{x,2},$ for rank three bundles and $x$ a degree one place. Let 
$$ g := \left( \begin{array}{ccc}
\pi_{x}^{d_1} &  &  \\ 
 & \pi_{x}^{d_2} &  \\ 
 &  & \pi_{x}^{d_3}
\end{array}\right),$$ 
which corresponds to  $\mathcal{O}_{\mathbb{P}^1}(d_1) \oplus \mathcal{O}_{\mathbb{P}^1}(d_2) \oplus \mathcal{O}_{\mathbb{P}^1}(d_3), \text{ with } d_1 \geq d_2 \geq d_3.$ 

For $w= [1:a:b] \in \mathrm{Gr}(1,3)(\kappa(x)) =\mathbb{P}^2(\mathbb{F}_q), a \neq 0$ we have
$$g \xi_w = \left( \begin{small}\begin{array}{ccc}
\pi_{x}^{d_1} &  &  \\ 
 & \pi_{x}^{d_2} &  \\ 
 &  & \pi_{x}^{d_3}
\end{array}\end{small}\right) \left(  \begin{small} \begin{array}{ccc}
\pi_{x}^{} &  & a \\ 
 & \pi_{x}^{} & b \\ 
 &  & 1
\end{array} \end{small}\right) = 
\left( \begin{small} \begin{array}{ccc}
\pi_{x}^{d_1+1} &  & a\pi_{x}^{d_1} \\ 
 & \pi_{x}^{d_2+1} & b \pi_{x}^{d_2} \\ 
 &  & \pi_{x}^{d_3}
\end{array} \end{small}\right)  $$ 
$$ \sim \left( \begin{small} \begin{array}{ccc}
 &  & 1 \\ 
 & 1 &   \\ 
1 &  & 
\end{array} \end{small}\right) \left( \begin{small} \begin{array}{ccc}
\pi_{x}^{d_1+1} &  & a\pi_{x}^{d_2} \\ 
 & \pi_{x}^{d_2+1} & b \pi_{x}^{d_2} \\ 
 &  & \pi_{x}^{d_3}
\end{array} \end{small}\right) =
\left(  \begin{small} \begin{array}{ccc}
 &  & \pi_{x}^{d_3} \\ 
 & \pi_{x}^{d_2+1} & b \pi_{x}^{d_2} \\ 
\pi_{x}^{d_1+1} &  & a\pi_{x}^{d_1}
\end{array}  \end{small}\right) $$  
$$\sim
\left( \begin{small} \begin{array}{ccc}
 &  & \pi_{x}^{d_3} \\ 
 & \pi_{x}^{d_2+1} & b \pi_{x}^{d_2} \\ 
\pi_{x}^{d_1+1} &  & a\pi_{x}^{d_1}
\end{array} \end{small} \right)  
\left( \begin{small} \begin{array}{ccc}
-a &  &  \\ 
 & 1 &  \\ 
\pi_{x} &  & 1
\end{array} \end{small}\right) $$ 
$$ =
\left( \begin{small} \begin{array}{ccc}
\pi_{x}^{d_3+1} &  & \pi_{x}^{d_3} \\ 
b\pi_{x}^{d_2+1} & \pi_{x}^{d_2+1} & b \pi_{x}^{d_2} \\ 
 &  & a\pi_{x}^{d_1}
\end{array} \end{small}\right) \sim
\left(  \begin{small}\begin{array}{ccc}
1 &  & -a^{-1}\pi_{x}^{d_3-d_1} \\ 
 & 1 &  \\ 
&  & 1
\end{array} \end{small}\right)
\left(  \begin{small} \begin{array}{ccc}
\pi_{x}^{d_3+1} &  & \pi_{x}^{d_3} \\ 
b\pi_{x}^{d_2+1} & \pi_{x}^{d_2+1} & b \pi_{x}^{d_2} \\ 
 &  & a\pi_{x}^{d_1}
\end{array} \end{small} \right) $$
$$ =
\left(  \begin{small} \begin{array}{ccc}
\pi_{x}^{d_3+1} &  &  \\ 
b\pi_{x}^{d_2+1} & \pi_{x}^{d_2+1} & b \pi_{x}^{d_2} \\ 
 &  & a\pi_{x}^{d_1}
\end{array} \end{small}\right) \sim 
\left(  \begin{small}\begin{array}{ccc}
1 &  &  \\ 
& 1 & -a^{-1}b \pi_{x}^{d_2-d_1} \\ 
 &  & 1
\end{array} \end{small}\right) 
\left( \begin{small} \begin{array}{ccc}
\pi_{x}^{d_3+1} &  &  \\ 
b\pi_{x}^{d_2+1} & \pi_{x}^{d_2+1} & b \pi_{x}^{d_2} \\ 
 &  & a\pi_{x}^{d_1}
\end{array} \end{small}\right)  $$ 
$$=\left( \begin{small} \begin{array}{ccc}
\pi_{x}^{d_3+1} &  &  \\ 
b\pi_{x}^{d_2+1} & \pi_{x}^{d_2+1} &  \\ 
 &  & a\pi_{x}^{d_1}
\end{array} \end{small}\right) \sim
\left(  \begin{small}\begin{array}{ccc}
\pi_{x}^{d_3+1} &  &  \\ 
b\pi_{x}^{d_2+1} & \pi_{x}^{d_2+1} &  \\ 
 &  & a\pi_{x}^{d_1}
\end{array} \end{small}\right)
\left(  \begin{small} \begin{array}{ccc}
1 &  &  \\ 
-b & 1 &  \\ 
 &  & a^{-1}
\end{array} \end{small} \right)  $$ 
$$=
\left( \begin{small} \begin{array}{ccc}
\pi_{x}^{d_3+1} &  &  \\ 
 & \pi_{x}^{d_2+1} &  \\ 
 &  & a\pi_{x}^{d_1}
\end{array} \end{small}\right) \sim 
\left(  \begin{small} \begin{array}{ccc}
\pi_{x}^{d_1} &  &  \\ 
 & \pi_{x}^{d_2+1} &  \\ 
 &  & \pi_{x}^{d_3+1}
\end{array} \end{small}\right) =: g_1 $$
with multiplicity $q(q-1).$ For $a=0$ and $b \neq 0$
$$g \xi_w = \left( \begin{small} \begin{array}{ccc}
\pi_{x}^{d_1+1} &  & \\ 
 & \pi_{x}^{d_2+1} & b \pi_{x}^{d_2} \\ 
 &  & \pi_{x}^{d_3}
\end{array} \end{small}\right)  \sim
\left( \begin{small} \begin{array}{ccc}
 &  & a\pi_{x}^{d_1+1} \\ 
b \pi_{x}^{d_2} & \pi_{x}^{d_2+1} &  \\ 
\pi_{x}^{d_3} &  &
\end{array} \end{small}\right) $$ 
$$ \sim 
\left(  \begin{small} \begin{array}{ccc}
 &  & a\pi_{x}^{d_1+1} \\ 
b \pi_{x}^{d_2} & \pi_{x}^{d_2+1} &  \\ 
\pi_{x}^{d_3} &  &
\end{array} \end{small}\right) 
\left(  \begin{small}\begin{array}{ccc}
1 & \pi_{x}^{} &  \\ 
 & -b &  \\ 
 &  & 1
\end{array} \end{small}\right) = 
\left(  \begin{small} \begin{array}{ccc}
 &  & a\pi_{x}^{d_1+1} \\ 
b \pi_{x}^{d_2} &  &  \\ 
\pi_{x}^{d_3} & \pi_{x}^{d_3+1}  &
\end{array} \end{small}\right)  $$ 
$$ \sim
\left( \begin{small} \begin{array}{ccc}
1 &  &  \\ 
 & 1 &  \\ 
 & -b^{-1} \pi_{x}^{d_3-d_2}  & 1
 \end{array} \end{small}\right) 
 \left(  \begin{small} \begin{array}{ccc}
 &  & a\pi_{x}^{d_1+1} \\ 
b \pi_{x}^{d_2} &  &  \\ 
\pi_{x}^{d_3} & \pi_{x}^{d_3+1} &
\end{array} \end{small}\right) =
\left( \begin{small} \begin{array}{ccc}
 &  & \pi_{x}^{d_1+1} \\ 
b \pi_{x}^{d_2} & &  \\ 
 & \pi_{x}^{d_3+1} &
\end{array} \end{small}\right) $$ 
$$ \sim 
 \left(  \begin{small} \begin{array}{ccc}
\pi_{x}^{d_1+1} &  &  \\ 
 & b \pi_{x}^{d_2} &  \\ 
 & & \pi_{x}^{d_3+1}
\end{array} \end{small}\right)   
 \left(  \begin{small}\begin{array}{ccc}
1 &  &  \\ 
 & b^{-1}  &  \\ 
 & & 1
\end{array} \end{small}\right) = 
  \left(  \begin{small} \begin{array}{ccc}
\pi_{x}^{d_1+1} &  &  \\ 
 &  \pi_{x}^{d_2} &  \\ 
 & & \pi_{x}^{d_3+1}
\end{array} \end{small} \right) =: g_2 
$$    
with multiplicity $q-1.$ If $a=b=0$, 
$$g \xi_w =  \left(  \begin{small} \begin{array}{ccc}
\pi_{x}^{d_1+1} &  &  \\ 
 & b \pi_{x}^{d_2+1} &  \\ 
 & & \pi_{x}^{d_3}
\end{array} \end{small}\right) =: g_3 $$
with multiplicity $1$. 

For $w= [0:0:1] \in \mathbb{P}^2(\mathbb{F}_q)$
$$g \xi_w =  \left( \begin{small} \begin{array}{ccc}
\pi_{x}^{d_1} &  &  \\ 
 & b \pi_{x}^{d_2} &  \\ 
 & & \pi_{x}^{d_3}
\end{array} \end{small}\right) 
 \left( \begin{small} \begin{array}{ccc}
1 &  &  \\ 
 & b \pi_{x} &  \\ 
 & & \pi_{x}
\end{array} \end{small} \right) = 
 \left( \begin{small} \begin{array}{ccc}
\pi_{x}^{d_1} &  &  \\ 
 & b \pi_{x}^{d_2+1} &  \\ 
 & & \pi_{x}^{d_3+1}
\end{array} \end{small}\right) = g_1   $$
with multiplicity $1.$ 

If $w=[0:1:a] \in \mathbb{P}^2(\mathbb{F}_q)$ and $a \neq 0$,
$$g \xi_w = \left( \begin{small} \begin{array}{ccc}
\pi_{x}^{d_1} &  &  \\ 
 & b \pi_{x}^{d_2} &  \\ 
 & & \pi_{x}^{d_3}
\end{array} \end{small}\right) \left( \begin{small} \begin{array}{ccc}
\pi_{x} & a &  \\ 
 & 1 &  \\ 
 & & \pi_{x}
\end{array} \end{small}\right) = 
\left(  \begin{small}\begin{array}{ccc}
\pi_{x}^{d_1+1} & a \pi_{x}^{d_1} &  \\ 
 & b \pi_{x}^{d_2} &  \\ 
 & & \pi_{x}^{d_3+1}
\end{array} \end{small}\right)  $$ 
$$\sim
\left(  \begin{small} \begin{array}{ccc}
\pi_{x}^{d_1+1} & a \pi_{x}^{d_1} &  \\ 
 & b \pi_{x}^{d_2} &  \\ 
 & & \pi_{x}^{d_3+1}
\end{array} \end{small} \right) 
\left(  \begin{small} \begin{array}{ccc}
-a &  &  \\ 
\pi_{x}^{} & 1 &  \\ 
 & & 1
\end{array} \end{small}\right) =
\left( \begin{small} \begin{array}{ccc}
 & a \pi_{x}^{d_1} &  \\ 
\pi_{x}^{d_2+1} & \pi_{x}^{d_2} &  \\ 
 & & \pi_{x}^{d_3+1}
\end{array} \end{small}\right)  $$ 
$$\sim
\left( \begin{small} \begin{array}{ccc}
1 & -a^{-1} \pi_{x}^{d_2-d_1} &  \\ 
 & 1 &  \\ 
 & & 1
\end{array} \end{small}\right)
\left( \begin{small} \begin{array}{ccc}
\pi_{x}^{d_2+1} & \pi_{x}^{d_2} &  \\ 
 & a \pi_{x}^{d_1} &  \\ 
 & & \pi_{x}^{d_3+1}
\end{array} \end{small}\right) = 
\left( \begin{small} \begin{array}{ccc}
\pi_{x}^{d_2+1} &  &  \\ 
 & a \pi_{x}^{d_1} &  \\ 
 & & \pi_{x}^{d_3+1}
\end{array} \end{small}\right)  $$ 
$$\sim
\left(  \begin{small}\begin{array}{ccc}
\pi_{x}^{d_1} &  &  \\ 
 &  \pi_{x}^{d_2+1} &  \\ 
 & & \pi_{x}^{d_3+1}
\end{array} \end{small}\right) = g_1$$
with multiplicity $q-1.$ For $a=0$,
$$g \xi_w = \left(  \begin{small}\begin{array}{ccc}
\pi_{x}^{d_1} &  &  \\ 
 &  \pi_{x}^{d_2} &  \\ 
 & & \pi_{x}^{d_3}
\end{array} \end{small}\right) \left( \begin{small} \begin{array}{ccc}
\pi_{x} &  &  \\ 
 & 1 &  \\ 
 & & \pi_{x}
\end{array} \end{small}\right)= \left(  \begin{small} \begin{array}{ccc}
\pi_{x}^{d_1+1} &  &  \\ 
 & b \pi_{x}^{d_2} &  \\ 
 & & \pi_{x}^{d_3+1}
\end{array} \end{small} \right) = g_2$$
with multiplicity $1.$ Therefore the graph of $\Phi_{x,2}$ is given by:

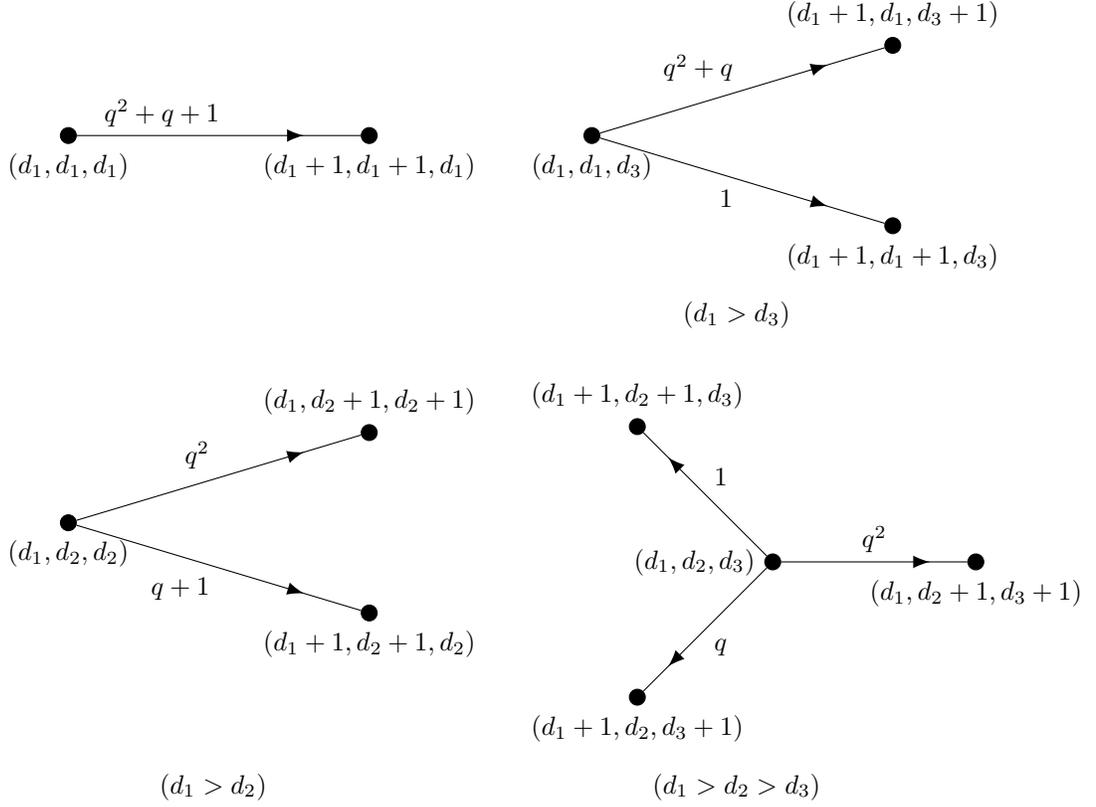
\begin{center}

\begin{figure}[h]
\centering

\begin{minipage}[b]{0.45\linewidth}

\[
  \beginpgfgraphicnamed{tikz/fig6}
  \begin{tikzpicture}[>=latex, scale=2]
        \vertex[circle,fill,label={below:${\footnotesize (d_1,d_1,d_1)}$}](00) at (0,0) {};
        \vertex[circle,fill,label={below:${\footnotesize (d_1+1,d_1+1,d_1)}$}](11) at (2,0) {};
           \path[-,font=\scriptsize]
   (00) edge[->-=0.8] node[pos=0.3,auto,black] {\normalsize  $q^2+q+1$} (11)
    ;
   \end{tikzpicture}
 \endpgfgraphicnamed
\]
\vspace*{1.2cm}
\end{minipage} 
\hfill
\begin{minipage}[b]{0.45\linewidth}
\[
  \beginpgfgraphicnamed{tikz/fig7}
  \begin{tikzpicture}[>=latex, scale=2]
        \vertex[circle,fill,label={below:$(d_1,d_1,d_3)$}](00) at (0,0) {};
        \vertex[circle,fill,label={above:$(d_1+1,d_1,d_3+1)$}](11) at (2,0.6) {};
        \vertex[circle,fill,label={below:$(d_1+1,d_1+1,d_3)$}](10) at (2,-0.6) {};
  
            \path[-,font=\scriptsize]
    (00) edge[->-=0.8] node[pos=0.5,auto,black,swap] {\normalsize $1$} (10)
    (00) edge[->-=0.8] node[pos=0.5,auto,black] {\normalsize  $q^2 + q$} (11)
    ;
   \end{tikzpicture}
 \endpgfgraphicnamed
\]
\begin{center}  $(d_1> d_3)$ \end{center}
\end{minipage}
\end{figure}

\end{center}
\begin{center}

\begin{figure}[h]
\centering
\begin{minipage}[b]{0.45\linewidth}
\[
  \beginpgfgraphicnamed{tikz/fig8}
  \begin{tikzpicture}[>=latex, scale=2]
        \vertex[circle,fill,label={below:$(d_1,d_2,d_2)$}](00) at (0,0) {};
        \vertex[circle,fill,label={above:$(d_1,d_2+1,d_2+1)$}](11) at (2,0.6) {};
        \vertex[circle,fill,label={below:$(d_1+1,d_2+1,d_2)$}](10) at (2,-0.6) {};
  
            \path[-,font=\scriptsize]
    (00) edge[->-=0.8] node[pos=0.5,auto,black,swap] {\normalsize $q+1$} (10)
    (00) edge[->-=0.8] node[pos=0.5,auto,black] {\normalsize  $q^2$} (11)
    ;
   \end{tikzpicture}
 \endpgfgraphicnamed
\]
\vspace{0.7cm}
\begin{center}  $(d_1> d_2)$ \end{center}
\end{minipage} 
\hfill
\begin{minipage}[b]{0.45\linewidth}
\[
  \beginpgfgraphicnamed{tikz/fig9}
  \begin{tikzpicture}[>=latex, scale=1.8]
        \vertex[circle,fill,label={left:${\tiny(d_1,d_2,d_3)}$}](00) at (0,0) {};
        \vertex[circle,fill,label={below:$(d_1,d_2+1,d_3+1)$}](11) at (1.5,0) {};
         \vertex[circle,fill,label={below:$(d_1+1,d_2,d_3+1)$}](12) at (-1,-1) {};
          \vertex[circle,fill,label={above:$(d_1+1,d_2+1,d_3)$}](13) at (-1,1) {};
        
           \path[-,font=\scriptsize]
   (00) edge[->-=0.8] node[pos=0.5,auto,black] {\normalsize  $q^2$} (11)
   (00) edge[->-=0.8] node[pos=0.5,auto,black] {\normalsize  $q$} (12)
   (00) edge[->-=0.8] node[pos=0.5,auto,black,swap] {\normalsize  $1$} (13)
    ;
   \end{tikzpicture}
 \endpgfgraphicnamed
\]
\begin{center}  $(d_1> d_2> d_3)$ \end{center}
\end{minipage}
\end{figure}

\end{center}


\end{example}


We proceed with the investigation of places of larger degree. Let us fix a place $x$ of degree one and let $y$ be a place of degree $d \geq 1.$ For determining the edges of $\mathcal{G}_{y,r}$,  we have to find the standard representative of $g \xi_w$ in $G(R^{-}) \setminus G(R^{\pm}) / G(R^{+})$ for $g= \pi_{x}^{\mathrm{\textbf{d}}} , 
\mathrm{\textbf{d}} \in  (\Z^n)_{+}$  and for every $w \in \mathrm{Gr}(n-r,n)(\kappa(y))$ where $\xi_w$ is as in Theorem \ref{maintheorem}.

The problem here is that $\xi_w$ has nontrivial entries in a place different from $x$ and, \textit{a priori}  we cannot use the reduction to a standard representative in $G(R^{-}) \setminus G(R^{\pm}) / G(R^{+}).$ Thus we have to find an equivalence class for $ \xi_w$  which depends only on the $x$-component. 

Let $S$ be the set of $(n \times n)$-matrices $(a_{ij})_{n\times n}$ defined as follows. Given $\lambda = (j_1, \ldots, j_{n-r}) \in J(n-r,n)$,  $(a_{ij})_{n \times n}$ is  an upper triangular matrix defined as follows:
\begin{itemize}
\item $a_{ii} = \pi_{x}^{d} \text{ if  $i \neq j_k$ for $k=1,\ldots,n-r$}   \text{ and }  a_{ii}=1 \text{ if $i$ is equals to some $j_k$}$  \\
\item $ a_{ij} := \sum_{l=0}^{d-1} a_{ij}^{l} \pi_{x}^{l} \in \mathbb{F}_q[\pi_x] \text{ if } i < j, a_{jj} = 1 \text{ and } a_{ii} = \pi_{x}^{d} $ \\
\item $a_{ij} = 0 \text{ otherwise.}$
\end{itemize}

Observe that we have $ \#\kappa(y) = q^d$ possibilities for the sum $ \sum_{l=0}^{d-1} a_{ij}^{l} \pi_{x}^{l}$ with $a_{ij}^{l} \in \mathbb{F}_q.$    

\begin{proposition}  \label{propforhightdegree}Keeping the above notation, the $\Phi_{y,r}-$neighbors of $[g]$ are the classes $[g \delta]$ with $\delta \in S.$ 

\begin{proof} We have to show that there is a bijection
$$\Psi : \{\xi_w | w \in \mathrm{Gr}(n-r,n)(\kappa(y))\} \stackrel{1:1}{\longrightarrow} S$$ 
such that $[g \xi_w] = [g \Psi(\xi_w)]$ in $G(F) \setminus G(\mathbb{A}) / K.$ 
First of all, we can suppose $\pi_y$ has a nontrivial valuation only in $y$ and $x.$ Indeed, let 
$$D(\pi_y) := \sum_{z \in |\mathbb{P}^1|} v_z(\pi_y) z$$
be the divisor associated with $\pi_y$, write 
$$D(\pi_y) = D + D' \text{ with } D := y-dx \text{ and } D' := d x + \sum_{z \neq y} v_z(\pi_y) z.$$ 
Since, $\deg D' =0$ and $|Cl^0 (F) |= 1$, we have $D(\pi_y) - D = D(f)$ for some $f \in F.$ Moreover $v_y(f)=0$ since $D(f) = D'$. Thus $f \in \mathcal{O}_{F,y}\setminus \{0\}$ and $\pi_y f^{-1}$ is a uniformizer for $y$ with
$$D(\pi_y f^{-1}) = D(\pi_y) - D(f) = D + D' - D' = y - dx.$$
Replacing $\pi_y$ by $\pi_y f^{-1}$, we can assume that $\pi_y$ has a nontrivial valuation only in $y$ and $x.$

Let $\delta \in G(F)$ denote the inverse of $(\xi_w)_y =: h.$ For all places $z \neq x,y,$ the canonical embedding
$$G(F) \hookrightarrow G(F_z)$$
sends $\delta$ to a matrix $\delta_z \in K_z$ since $v_z(\pi_y) =0.$ Let $k \in K$ such that $(k)_z = \delta_z$ for $z \neq x,y$, $(k)_x = I_n$ and $(k)_y = I_n$. Let $\delta_w = \delta \xi_w k^{-1} \in G(\mathbb{A}).$ Observe that  only the $x-$component of $\delta_w$ is nontrivial, with  $(\delta_w)_x= h^{-1}.$ We have
$$[g \xi_w] = [\delta g \xi_w k^{-1}] = [g \delta \xi_w k^{-1}] = [g \delta_x]$$
where the first equality holds because $\delta \in G(F)$ and $k^{-1} \in K$, the second because $g$ and $\delta$ are diagonal matrices and the last due to the definition of $\delta_w.$ 

To finish the proof, we need to show that $\delta_w \in S$ and $S= \{\delta_w \}$. The only nontrivial component of $\delta_w$ is $(\delta_w)_x = (\xi_w)_{y}^{-1}.$ Since $v_x(\pi_y) = -d$, we have $v_x(\pi_y) =d$ and $\pi_{y}^{-1} = u \pi_{x}^{d}$ for some $u \in \mathcal{O}_{F,x}\setminus \{0\}.$ Thus the diagonal of $(\delta_w)_x$ has entries $u \pi_{x}^{d}$ at the positions that $(\xi_w)_y$ has entries $\pi_y$.  Using the reduction from Lemma \ref{hellslemma}, we have $(\delta_w)_z = I_n$ for $z \neq x$ and conclude that $(\delta_w)_x$ is as desired. In the end, since each $w \in \mathrm{Gr}(n-r,n)(\kappa(y))$ gives us a unique $\xi_w$ and each $\xi_w$ is associated with some $\delta_w \in S$ and  $\#S = \#\{\xi_w\}$, $\Psi$ is a bijection. Each $\delta_w$ is in a different class, i.e.\ the matrices as defined before occur as the $x-$components of an associated $\delta_w.$
\end{proof}
\end{proposition}

By the previous proposition, the graph $\mathcal{G}_{y,r}$ depends only on the degree of $y$. Since we have a representation for $\xi_w$ as $\delta_w,$ whose only nontrivial component is the $x$-component, we can use the same reduction as before to find the neighbors of a vertex in the graph of $\Phi_{y,r}$. Note that if $y$ is a place of degree one $x$, then $\mathcal{G}_{y,r}$ is the same graph as $\mathcal{G}_{x,r}.$ 

\begin{example} Let us calculate the graph $\mathcal{G}_{y,1}$ for rank $2$-bundles and $\deg(y)=2$. According the last proposition, we have to find the standard representatives for 
$$\left( \begin{small}
\begin{array}{cc}
\pi_{x}^{d_1} &  \\ 
 & \pi_{x}^{d_2}
\end{array} \end{small} \right) \left( \begin{small}\begin{array}{cc}
1 &  \\ 
 & \pi_{x}^{2}
\end{array}\end{small} \right) \text{ and } \left(\begin{small} \begin{array}{cc}
\pi_{x}^{d_1} &  \\ 
 & \pi_{x}^{d_2}
\end{array} \end{small} \right)\left(\begin{small} \begin{array}{cc} \pi_{x}^{2} & a_0 + a_1 \pi_{x}^{} \\ & 1 \end{array}\end{small} \right)$$ 
where $d_1 \geq d_2$ and $a_0, a_1 \in \mathbb{F}_q.$

In the first case, we have
$$\left( \begin{array}{cc}
\pi_{x}^{d_1} &  \\ 
 & \pi_{x}^{d_2}
\end{array} \right) \left( \begin{array}{cc}
1 &  \\ 
 & \pi_{x}^{2}
\end{array} \right) = \left( \begin{array}{cc}
\pi_{x}^{d_1} &  \\ 
 & \pi_{x}^{d_2+2}
\end{array} \right)$$
with multiplicity one. In the second case, we first assume $a_0 \neq 0$ and let $s=a_0 + a_1 \pi_{x}^{} \in \mathcal{O}_{F,x}^{*}.$ Then
$$\left(\begin{small} \begin{array}{cc}
\pi_{x}^{d_1} &  \\ 
 & \pi_{x}^{d_2}
\end{array} \end{small}\right)\left(\begin{small} \begin{array}{cc} \pi_{x}^{2} & a_0 + a_1 \pi_{x}^{} \\ & 1 \end{array}\end{small} \right) =
\left( \begin{small} \begin{array}{cc}
\pi_{x}^{d_1+2} & s \pi_{x}^{d_1} \\ 
 & \pi_{x}^{d_2}
\end{array} \end{small}\right)  $$  
$$ \sim \left( \begin{small} \begin{array}{cc}
1 &  -s^{-1} \pi_{x}^{d_2-d_1} \\ 
 & 1
\end{array}\end{small} \right) \left(\begin{small} \begin{array}{cc}
\pi_{x}^{d_2} &  \\ 
 s \pi_{x}^{d_1}& \pi_{x}^{d_1+2}
\end{array}\end{small} \right) = \left( \begin{small} \begin{array}{cc}
 & -s^{-1} \pi_{x}^{d_2+2} \\ 
 s \pi_{x}^{d_1} & \pi_{x}^{d_1+2}
\end{array}\end{small} \right) $$ 
$$ \sim  \left(\begin{small} \begin{array}{cc}
 s \pi_{x}^{d_1} & \pi_{x}^{d_1+2}\\ 
  & -s^{-1} \pi_{x}^{d_2+2}
\end{array} \end{small}\right)\left( \begin{small} \begin{array}{cc}
1 & -s^{-1}\pi_{x}^{2} \\ 
 & 1
\end{array}\end{small} \right) = \left( \begin{small} \begin{array}{cc}
s \pi_{x}^{d_1} &  \\ 
 & -s^{-1}\pi_{x}^{d_2+2}
\end{array} \end{small} \right) \sim \left(\begin{small}  \begin{array}{cc}
\pi_{x}^{d_1} &  \\ 
 & \pi_{x}^{d_2+2}
\end{array} \end{small} \right)$$
whose multiplicity is $q(q-1).$ If $a_0 =0$ and $a_1 \neq 0,$ then 
$$\left( \begin{small} \begin{array}{cc}
\pi_{x}^{d_1} &  \\ 
 & \pi_{x}^{d_2}
\end{array}\end{small} \right)\left( \begin{small} \begin{array}{cc} \pi_{x}^{2} & a_0 + a_1 \pi_{x}^{} \\ & 1 \end{array} \end{small}\right) =
\left( \begin{small} \begin{array}{cc}
\pi_{x}^{d_1+2} & a_1 \pi_{x}^{d_1+1} \\ 
 & \pi_{x}^{d_2}
\end{array}\end{small} \right)  $$
$$\sim \left( \begin{small} \begin{array}{cc}
1 & -a_{1}^{-1} \pi_{x}^{d_2-d_1-1} \\ 
 & 1
\end{array} \end{small}\right) \left( \begin{small} \begin{array}{cc}
\pi_{x}^{d_2} &  \\ 
a_1 \pi_{x}^{d_1+1} & \pi_{x}^{d_1+2}
\end{array}\end{small} \right) = \left( \begin{small}\begin{array}{cc}
0 & -a_{1}^{-1} \pi_{x}^{d_2+1} \\ 
a_1 \pi_{x}^{d_1+1} & \pi_{x}^{d_1+2}
\end{array}\end{small} \right) $$ 
$$\sim  \left(\begin{small} \begin{array}{cc}
a_1\pi_{x}^{d_1+1} & \pi_{x}^{d_1+2} \\ 
 & -a_{1}^{-1}\pi_{x}^{d_2+1}
\end{array} \end{small}\right) \left( \begin{small}\begin{array}{cc}
1 & -a_{1}^{-1} \pi_{x}^{} \\ 
 & 1
\end{array} \end{small}\right) $$ 
$$= \left(\begin{small} \begin{array}{cc}
a_1\pi_{x}^{d_1+1} &  \\ 
 & -a_{1}^{-1} \pi_{x}^{d_2+1}
\end{array}\end{small} \right) \sim \left( \begin{small} \begin{array}{cc}
\pi_{x}^{d_1+1} &  \\ 
 & \pi_{x}^{d_2+1}
\end{array} \end{small}\right)$$
with multiplicity $q-1.$ If $a_1 = b_1 =0$, then 
$$\left( \begin{small}\begin{array}{cc}
\pi_{x}^{d_1} &  \\ 
 & \pi_{x}^{d_2}
\end{array} \end{small}\right)\left( \begin{small}\begin{array}{cc} \pi_{x}^{2} & a_0 + a_1 \pi_{x}^{} \\ & 1 \end{array} \end{small}\right) =
\left( \begin{small}\begin{array}{cc}
\pi_{x}^{d_1+2} &  \\ 
 & \pi_{x}^{d_2}
\end{array} \end{small}\right)  $$ 
with multiplicity $1$. Therefore, the graph $\mathcal{G}_{y,1}$ is as follows.

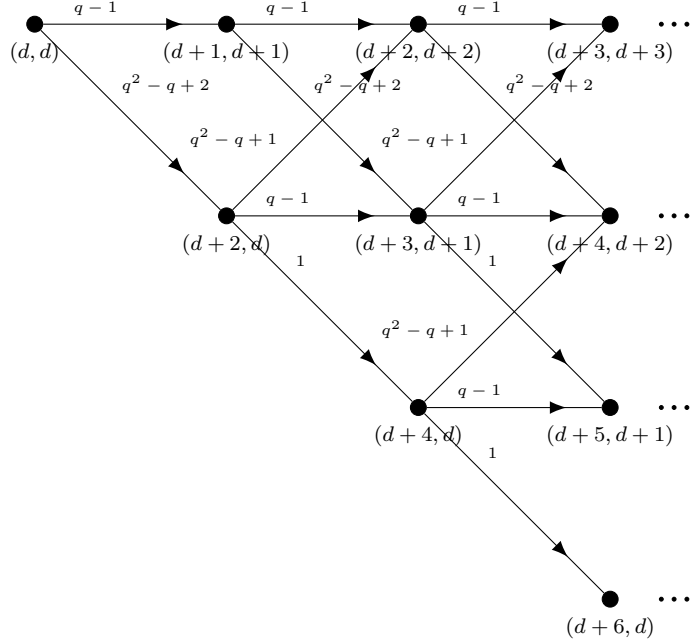
\begin{figure}[h]
  \beginpgfgraphicnamed{tikz/graphd2}
  \begin{tikzpicture}[>=latex, scale=1.7]
        \vertex[circle,fill,label={below: {\footnotesize $(d,d)$}}](00) at (-1.5,0) {};
        \vertex[circle,fill,label={below: {\footnotesize $(d+1,d+1)$}}](11) at (0,0) {};
        \vertex[circle,fill,label={below: {\footnotesize $(d+2,d+2)$}}](10) at (1.5,0) {};
        \vertex[circle,fill,label={below: {\footnotesize $(d+3,d+3)$}}](22) at (3,0) {};
       
        \vertex[circle,fill,label={below: {\footnotesize $(d+2,d)$}}](20) at (0,-1.5) {};
        \vertex[circle,fill,label={below: {\footnotesize $(d+3,d+1)$}}](21) at (1.5,-1.5) {};
         \vertex[circle,fill,label={below: {\footnotesize $(d+4,d+2)$}}](23) at (3,-1.5) {};

        \vertex[circle,fill,label={below:{\footnotesize $(d+4,d)$}}](31) at (1.5,-3) {};
        \vertex[circle,fill,label={below:{\footnotesize $(d+5,d+1)$}}](32) at (3,-3) {};

        \vertex[circle,fill,label={below:{\footnotesize $(d+6,d)$}}](40) at (3,-4.5) {};
    \path[-,font=\scriptsize]
    (00) edge[->-=0.8] node[pos=0.3,auto,black] {\tiny $q-1$} (11)
    (11) edge[->-=0.8] node[pos=0.3,auto,black] {\tiny $q-1$} (10)
    (10) edge[->-=0.8] node[pos=0.3,auto,black] {\tiny $q-1$} (22)
   
    (20) edge[->-=0.8] node[pos=0.3,auto,black] {\tiny $q-1$} (21)
    (21) edge[->-=0.8] node[pos=0.3,auto,black] {\tiny $q-1$} (23)
    
    (31) edge[->-=0.8] node[pos=0.3,auto,black] {\tiny $q-1$} (32)
    
    (00) edge[->-=0.8] node[pos=0.4,auto,black] {\tiny $q^2-q+2$} (20)
     (11) edge[->-=0.8] node[pos=0.4,auto,black] {\tiny $q^2-q+2$} (21)
      (10) edge[->-=0.8] node[pos=0.4,auto,black] {\tiny $q^2-q+2$} (23)
     
       (20) edge[->-=0.8] node[pos=0.3,auto,black] {\tiny $q^2-q+1$} (10)
       (21) edge[->-=0.8] node[pos=0.3,auto,black] {\tiny $q^2-q+1$} (22)     
       (21) edge[->-=0.8] node[pos=0.3,auto,black] {\tiny $1$} (32)  
       (20) edge[->-=0.8] node[pos=0.3,auto,black] {\tiny $1$} (31)

(31) edge[->-=0.8] node[pos=0.3,auto,black] {\tiny $q^2-q+1$} (23)       
  (31) edge[->-=0.8] node[pos=0.3,auto,black] {\tiny $1$} (40)        
       ;
        \draw (3.4,0) circle (0.015cm);
\fill  (3.4,0) circle (0.015cm);
  \draw (3.5,0) circle (0.015cm);
\fill  (3.5,0) circle (0.015cm);  
\draw (3.6,0) circle (0.015cm);
\fill  (3.6,0) circle (0.015cm);

        \draw (3.4,-1.5) circle (0.015cm);
\fill  (3.4,-1.5) circle (0.015cm);
  \draw (3.5,-1.5) circle (0.015cm);
\fill  (3.5,-1.5) circle (0.015cm);  
\draw (3.6,-1.5) circle (0.015cm);
\fill  (3.6,-1.5) circle (0.015cm);

        \draw (3.4,-3) circle (0.015cm);
\fill  (3.4,-3) circle (0.015cm);
  \draw (3.5,-3) circle (0.015cm);
\fill  (3.5,-3) circle (0.015cm);  
\draw (3.6,-3) circle (0.015cm);
\fill  (3.6,-3) circle (0.015cm);

   \draw (3.4,-4.5) circle (0.015cm);
\fill  (3.4,-4.5) circle (0.015cm);
  \draw (3.5,-4.5) circle (0.015cm);
\fill  (3.5,-4.5) circle (0.015cm);  
\draw (3.6,-4.5) circle (0.015cm);
\fill  (3.6,-4.5) circle (0.015cm);
\end{tikzpicture}
 \endpgfgraphicnamed
\end{figure}

\end{example} 

\newpage

\begin{example} We conclude with the description of the graph $\mathcal{G}_{y,2}$, for rank $3$ bundles where $y$ is a place of degree two. We do not write out all calculations, which are similar to the ones for $\mathcal{G}_{y,1}$. By Proposition  \ref{propforhightdegree}, we have to find the standard representatives for $g \delta_{w}$ with 
$$g:= \left( \begin{small}\begin{array}{ccc}
\pi_{x}^{d_1} &  & \\ 
 & \pi_{x}^{d_2}  & \\
 & & \pi_{x}^{d_3}
\end{array} \end{small}\right), \text{ and $\delta_w$ one of the matrices}  $$
$$
\left( \begin{small}\begin{array}{ccc}
1 &  & \\ 
 & \pi_{x}^{2}  & \\
 & & \pi_{x}^{2}
\end{array} \end{small}\right), 
\left( \begin{small} \begin{array}{ccc}
\pi_{x}^{2} & a_0 + a_1 \pi_{x}^{} & \\ 
 & 1  & \\
 & & \pi_{x}^{2}
\end{array}\end{small} \right) \text{ or } 
\left(\begin{small} \begin{array}{ccc}
\pi_{x}^{2} &  & a_0 + a_1 \pi_{x}^{}\\ 
 & \pi_{x}^{2}  & b_0 + b_1 \pi_{x}^{} \\
 & & 1
\end{array}\end{small} \right),$$
where $d_1 \geq d_2 \geq d_3$ and $a_i,b_i \in \mathbb{F}_q$ for $ i=0,1.$ As before,  we will write $(d_1,d_2,d_3)$ for the matrix $g$. 
 
In the first case, we obtain 
$$g \delta_w = \left(\begin{small} \begin{array}{ccc}
\pi_{x}^{d_1} &  & \\ 
 & \pi_{x}^{d_2}  & \\
 & & \pi_{x}^{d_3}
\end{array} \end{small}\right)
\left(\begin{small} \begin{array}{ccc}
1 &  & \\ 
 & \pi_{x}^{2}  & \\
 & & \pi_{x}^{2}
\end{array} \end{small}\right) \sim (d_1,d_2+2, d_3+2)
$$
with multiplicity $1$. In the second case, we have
$$g \delta_w = \left(\begin{small} \begin{array}{ccc}
\pi_{x}^{d_1} &  & \\ 
 & \pi_{x}^{d_2}  & \\
 & & \pi_{x}^{d_3}
\end{array} \end{small}\right)
\left( \begin{small} \begin{array}{ccc}
\pi_{x}^{2} & a_0 + a_1 \pi_{x}^{} & \\ 
 & 1  & \\
 & & \pi_{x}^{2}
\end{array} \end{small}\right) = 
\left( \begin{small} \begin{array}{ccc}
\pi_{x}^{d_1+2} & (a_0 + a_1 \pi_{x}^{}) \pi_{x}^{d_1} & \\ 
 & \pi_{x}^{d_2}  & \\
 & & \pi_{x}^{d_3+2}
\end{array}\end{small} \right)$$
and we have to analyse the following subcases.

$\bullet \   \ a_0 \neq 0$: \text{then } $g \delta_w \sim (d_1,d_2 + 2d_3+2)$ with multiplicity $q(q-1).$

$\bullet \   \  a_0 =0$ and $a_1=0$: \text{ then } $g \delta_w \sim (d_1+2,d_2,d_3+2)$ with multiplicity $1$. 

$\bullet \   \  a_0 =0$ and $a_1 \neq 0:$ \text{ then }$g \delta_w \sim (d_1+1,d_2+1,d_3+2)$ with multiplicity $1$.

In the third case, we have
$$g \delta_w =  \left( \begin{small} \begin{array}{ccc}
\pi_{x}^{d_1} &  & \\ 
 & \pi_{x}^{d_2}  & \\
 & & \pi_{x}^{d_3}
\end{array} \end{small}\right) \left(\begin{small} \begin{array}{ccc}
\pi_{x}^{2} &  & a_0 + a_1 \pi_{x}^{}\\ 
 & \pi_{x}^{2}  & b_0 + b_1 \pi_{x}^{} \\
 & & 1
\end{array}\end{small} \right) = \left(\begin{small} \begin{array}{ccc}
\pi_{x}^{d_1+2} &  & (a_0 + a_1 \pi_{x}^{})\pi_{x}^{d_1}\\ 
 & \pi_{x}^{d_2+2}  & (b_0 + b_1 \pi_{x}^{})\pi_{x}^{d_2} \\
 & & \pi_{x}^{d_3}
\end{array} \end{small}\right) $$ 
and consider the following subcases:
\begin{itemize}

\item $ a_0 \neq 0$ and $b_0 \neq 0:$ \text{ then } $g \delta_w \sim (d_1,d_2+ 2,d_3+2)$ with multiplicity $q^2(q-1)^2$.

\item $  a_0 =0, a_1 \neq 0$ and $b_0 \neq 0:$ \text{ then } $g \delta_w \sim (d_1+1,d_2+1,d_3+2)$ with multiplicity $q(q-1)^2$.

\item $  a_0 =0, a_1 = 0$ and $b_0 \neq 0:$ \text{ then }  $g \delta_w \sim (d_1+2,d_2,d_3+2)$ with multiplicity $q(q-1)$.

\item $  a_0 =0, a_1 \neq 0$ and $b_0 = 0, b_1 \neq 0:$ \text{ then } $g \delta_w \sim (d_1+1,d_2+2,d_3+1)$ with multiplicity $(q-1)^2$.

\item $  a_0 =0, a_1 \neq 0$ and $b_0 = 0, b_1=0:$ \text{ then } $g \delta_w \sim (d_1+1,d_2+2,d_3+1)$ with multiplicity $(q-1)$.

\item $  a_0 =0, a_1 = 0$ and $b_0 = 0, b_1 \neq 0:$ \text{ then } $g \delta_w \sim (d_1+2,d_2+1,d_3+1)$ with multiplicity $(q-1)$.

\item $ a_0 =0, a_1 = 0$ and $b_0 = 0, b_1 =0:$ \text{ then } $g \delta_w \sim (d_1+2,d_2+2,d_3)$ with multiplicity $1$.

\item $ a_0 \neq 0, $ and $b_0 = 0, b_1 =0:$ \text{ then } $g \delta_w \sim (d_1,d_2+2,d_3+2)$ with multiplicity $q(q-1)$.

\item $  a_0 \neq 0, $ and $b_0 = 0, b_1 \neq 0:$ \text{ then } $g \delta_w \sim (d_1,d_2+2,d_3+2)$ with multiplicity $q(q-1)^2$.

\end{itemize}

The graph $\mathcal{G}_{y,2}$ can be illustrated as follows.

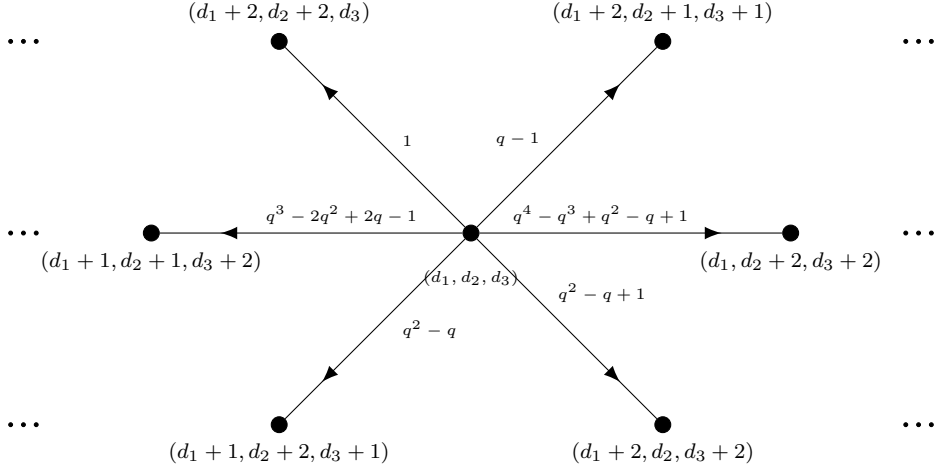
\begin{figure}[h]
  \beginpgfgraphicnamed{tikz/graph32}
  \begin{tikzpicture}[>=latex, scale=1.7]
        \vertex[circle,fill,label={[label distance=0.25cm]below: {\tiny $(d_1,d_2,d_3)$}}](00) at (0,0) {};

        \vertex[circle,fill,label={below: {\footnotesize $(d_1,d_2+2,d_3+2)$}}](10) at (2.5,0) {};
        \vertex[circle,fill,label={below: {\footnotesize $(d_1+1,d_2+1,d_3+2)$}}](11) at (-2.5,0) {};
        
        \vertex[circle,fill,label={below: {\footnotesize $(d_1+2,d_2,d_3+2)$}}](20) at (1.5,-1.5) {};
       \vertex[circle,fill,label={above: {\footnotesize $(d_1+2,d_2+2,d_3)$}}](21) at (-1.5,1.5) {};
        
        \vertex[circle,fill,label={below: {\footnotesize $(d_1+1,d_2+2,d_3+1)$}}](30) at (-1.5,-1.5) {};
         \vertex[circle,fill,label={above: {\footnotesize $(d_1+2,d_2+1,d_3+1)$}}](31) at (1.5,1.5) {};

    \path[-,font=\scriptsize]
    (00) edge[->-=0.8] node[pos=0.4,auto,black] {\tiny $q^4-q^3+q^2-q+1$} (10)
    (00) edge[->-=0.8] node[pos=0.4,auto,black,swap] {\tiny $q^3-2q^2+2q-1$} (11)

    (00) edge[->-=0.8] node[pos=0.4,auto,black] {\tiny $q^2-q+1$} (20)
    (00) edge[->-=0.8] node[pos=0.4,auto,black,swap] {\tiny $1$} (21)
    
    (00) edge[->-=0.8] node[pos=0.4,auto,black] {\tiny $q^2-q$} (30)
    (00) edge[->-=0.8] node[pos=0.4,auto,black] {\tiny $q-1$} (31)

       ;
        \draw (3.4,0) circle (0.015cm);
\fill  (3.4,0) circle (0.015cm);
  \draw (3.5,0) circle (0.015cm);
\fill  (3.5,0) circle (0.015cm);  
\draw (3.6,0) circle (0.015cm);
\fill  (3.6,0) circle (0.015cm);

        \draw (3.4,-1.5) circle (0.015cm);
\fill  (3.4,-1.5) circle (0.015cm);
  \draw (3.5,-1.5) circle (0.015cm);
\fill  (3.5,-1.5) circle (0.015cm);  
\draw (3.6,-1.5) circle (0.015cm);
\fill  (3.6,-1.5) circle (0.015cm);

        \draw (-3.4,0) circle (0.015cm);
\fill  (-3.4,0) circle (0.015cm);
  \draw (-3.5,0) circle (0.015cm);
\fill  (-3.5,0) circle (0.015cm);  
\draw (-3.6,0) circle (0.015cm);
\fill  (-3.6,0) circle (0.015cm);

        \draw (-3.4,-1.5) circle (0.015cm);
\fill  (-3.4,-1.5) circle (0.015cm);
  \draw (-3.5,-1.5) circle (0.015cm);
\fill  (-3.5,-1.5) circle (0.015cm);  
\draw (-3.6,-1.5) circle (0.015cm);
\fill  (-3.6,-1.5) circle (0.015cm);

        \draw (3.4,1.5) circle (0.015cm);
\fill  (3.4,1.5) circle (0.015cm);
  \draw (3.5,1.5) circle (0.015cm);
\fill  (3.5,1.5) circle (0.015cm);  
\draw (3.6,1.5) circle (0.015cm);
\fill  (3.6,1.5) circle (0.015cm);

        \draw (-3.4,1.5) circle (0.015cm);
\fill  (-3.4,1.5) circle (0.015cm);
  \draw (-3.5,1.5) circle (0.015cm);
\fill  (-3.5,1.5) circle (0.015cm);  
\draw (-3.6,1.5) circle (0.015cm);
\fill  (-3.6,1.5) circle (0.015cm);
   
\end{tikzpicture}
 \endpgfgraphicnamed
\end{figure}
Observe that, the sum up of multiplicities of edges originating in $[g]$ is $q^4 + q^2 + 1$ which is equal to $\# \mathrm{Gr}(1,3)(\mathbb{F}_{q^2}),$ as stated in Theorem \ref{maintheorem}.

\end{example}

\vspace*{1cm}

\textbf{Acknowledgements:} This article is part of the author's Ph.D thesis at
IMPA under the supervision of Oliver Lorscheid. He deeply
thanks him for his constant support, patience, encouragement and availability. He also would like to thank the referee for the carefully reading and several helpful comments. 
Funding: This work was supported by  Faperj [grant 200.322/2016];  and Fapesp [grant number 2017/21259-3].

\bibliographystyle{plain}

\vspace*{0.5cm}

\end{document}